\documentclass[12pt]{amsart}
\usepackage{amssymb}
\usepackage{graphicx} 
\usepackage{float}
\usepackage{xurl} 
\usepackage{mathrsfs}
\usepackage{xcolor}
\usepackage[all,color]{xy}
\usepackage{amsthm}
\usepackage{bbold}
\usepackage{amsmath}
\usepackage{a4wide}
\usepackage[
    colorlinks=true,
    hypertexnames=false,
    linkcolor=blue!50!black,
    citecolor=blue!50!black,
    urlcolor=blue!50!black
]{hyperref}

\newcommand{\leanrevision}{ca3135e88ade561a58f1b6f9ce57a0fd8cd96d2f}
\newcommand{\leanref}[3]{%
  \emergencystretch=1em\relax
  \href{https://github.com/GonzMA/FreeProductsVecLat/blob/\leanrevision/#1\#L#2}{%
    \textnormal{\footnotesize\ensuremath{\boldsymbol{\forall}}\,\nolinkurl{#3}}}}

\newcommand{\vertiii}[1]{{\left\vert\kern-0.25ex\left\vert\kern-0.25ex\left\vert #1
    \right\vert\kern-0.25ex\right\vert\kern-0.25ex\right\vert}}

\DeclareMathOperator*{\bigfree}{\raisebox{-0.6ex}{\scalebox{2}{$\ast$}}}

\newtheorem{theo}{Theorem}[section]
\newtheorem{lema}[theo]{Lemma}
\newtheorem{cor}[theo]{Corollary}
\newtheorem{prop}[theo]{Proposition}
\theoremstyle{definition}
\newtheorem{defin}[theo]{Definition}

\newtheorem{obs}[theo]{Observation}
\newtheorem{ex}[theo]{Example}
\DeclareMathOperator{\lat}{lat}
\DeclareMathOperator{\Van}{Van}

\title{Free Products of Vector and Banach Lattices}

\author[G. Mart\'inez-Fern\'andez]{Gonzalo Mart\'inez-Fern\'andez}
\address{Instituto de Ciencias Matem\'aticas\\
Consejo Superior de Investigaciones Cient\'ificas\\
C/ Nicol\'as Cabrera 13--15\\
28049 Madrid, Spain}
\email{gonzalo.martinez@icmat.es}

\author[P. Tradacete]{Pedro Tradacete}
\address{Instituto de Ciencias Matem\'aticas\\
Consejo Superior de Investigaciones Cient\'ificas\\
C/ Nicol\'as Cabrera 13--15\\
28049 Madrid, Spain}
\email{pedro.tradacete@icmat.es}

\date{}

\begin{document}

\begin{abstract}
We study free products of vector lattices and their relation to free
products of Banach lattices. We characterize the Archimedean vector
lattices whose vector lattice free product with every Archimedean vector lattice
is again Archimedean: these are precisely those whose finitely generated
sublattices are projective. For finitely generated lattices, this connects
the free-product property with polyhedral geometry through Baker--Beynon
duality. We also show that the Banach lattice free product is the completion of
the Archimedean vector lattice free product under the largest lattice norm
making the canonical maps contractive. For finitely many factors, the
Archimedean vector lattice free product is order dense in the Banach lattice
free product. Moreover, it is shown that bounded families of injective lattice homomorphisms
also induce injective homomorphisms between Banach lattice free products.
Finally, in analogy with the theory of separably injective
Banach lattices, we study relative injectivity and
universal disposition for the class of vector lattices with countably many
proper prime ideals. Every vector lattice embeds into one of universal
disposition for this class, while ordered ultrapowers of the real line
provide relatively injective examples of cardinality continuum.
\end{abstract}



\maketitle

\tableofcontents

\section{Introduction}
Free products, also called coproducts, describe how to combine a family of
objects in a way that allows morphisms defined on each of the individual factors to determine uniquely a morphism on the resulting object. In many familiar categories, these constructions include disjoint unions of sets and
topological spaces, direct sums of vector spaces and modules, or free
products of groups. In the category of Banach spaces and contractive linear
operators, the corresponding construction coincides with the $\ell_1$-sum. Although
these examples share the same universal property, the structure of a free
product depends fundamentally on the class of morphisms under consideration. For
vector lattices, when morphisms are the maps preserving the linear and lattice operations, the direct sum with coordinatewise order does not have
the required property: its distinct summands are disjoint, whereas the
images of independently chosen lattice homomorphisms need not be. A free
product must therefore allow lattice operations between elements from different factors without imposing disjointness or any other prescribed relations among them.

The study of these constructions has a substantial history in ordered
algebra and functional analysis. The following developments provide
background for the questions considered here: In vector lattice theory, Topping~\cite{Topping1965} studied projective and injective objects and proved that there are no nonzero
injective objects in the category of vector lattices and linear lattice
homomorphisms. This obstruction motivates extension problems relative to
restricted classes of embeddings (see Section~\ref{sec:countable-spectrum}). Free products of lattice-ordered groups were studied by Holland
and Scrimger~\cite{HollandScrimger1972} and, in varieties of
lattice-ordered groups and in the abelian case, by
Martinez~\cite{Martinez1972FreeProducts,Martinez1973AbelianFreeProducts}.
Powell and Tsinakis~\cite{PowellTsinakis1983} subsequently gave concrete
constructions in the abelian case and proved that embeddings of the factors induce embeddings of their free products.
This can be seen as the group-theoretic counterpart of
Corollary~\ref{cor:freeprod-vl-injective}. The amalgamation results of
Pierce~\cite{Pierce:1973,Pierce:1976} and the treatment of push-outs by
Hager and Madden~\cite{HagerMadden1983} connect these constructions with
extension properties of lattice homomorphisms. Broader accounts of free
products and amalgamation appear in the chapters by Powell and
Tsinakis~\cite{PowellTsinakis1989FreeProducts} and by Tsinakis and
Powell~\cite{TsinakisPowell1989Amalgamations}. See also \cite{Mundici1988,Ton1998,HagerKimber2014} for further developments.

In the Banach lattice setting, isometric push-outs were constructed in \cite{amalg} and used to obtain amalgamation and relative injectivity results.
Tursi~\cite{Tursi2023} connected amalgamation with metric Fra\"iss\'e
theory to construct a separable approximately ultrahomogeneous Banach
lattice that is isometrically universal for separable Banach lattices.
Also, Garc\'ia-S\'anchez, Leung, Taylor and
the second author~\cite{GarciaSanchezLeungTaylorTradacete2025UpperP} adapted
push-out constructions in classes of Banach lattices subject to
convexity conditions, including $p$-convex lattices and lattices with
upper $p$-estimates, together with quantitative results on preservation
of embeddings. This connects the categorical constructions with the
geometry of lattice norms.

In our previous paper \cite{freeprods_BL}, we constructed arbitrary free
products of Banach lattices where morphisms are lattice homomorphisms. These can be described as appropriate quotients of free Banach lattices and we could
establish several basic results such as their compatibility with isometric embeddings, quotients and
projectivity. For nonempty compact Hausdorff spaces $K,L$, we also obtained the representation
\[
 C(K)*_{\scriptscriptstyle{\mathcal{BL}}}C(L)\simeq C(K*L),
\]
where $K*L$ denotes the topological join, together with an explicit formula for
the free-product norm. Through the double suspension theorem, from geometric topology, this representation yields surprising examples of free factors of $\operatorname{FBL}[\ell_1]$ that are not themselves lattice isomorphic to free Banach lattices.

The purpose of the present paper is to understand the vector-lattice
structure underlying these constructions and its relation to Banach
lattice theory. Three categories enter naturally: vector lattices
$\mathcal{VL}$, Archimedean vector lattices $\mathcal{AVL}$, both with lattice
homomorphisms, and Banach lattices $\mathcal{BL}$ with contractive lattice
homomorphisms. The relations between these categories raise several questions. When is the
ordinary vector-lattice free product of Archimedean factors still
Archimedean? Which elements are identified on passing to the Banach lattice
free product? Do injective homomorphisms remain injective after taking free
products, even when their ranges are not closed? These questions concern
both the algebraic relations imposed by the order and the additional
structure introduced by a lattice norm and its completion.

We begin with the construction of free products and push-outs in
$\mathcal{VL}$. Prime ideals reduce the study of embeddings to totally
ordered quotients. This gives a direct proof of the classical fact that
push-outs preserve injective lattice homomorphisms, and hence that free
products preserve embeddings. The arguments use the theory of amalgamation of lattice-ordered groups and vector lattices developed in
\cite{Pierce:1973,Pierce:1976,HagerMadden1983}. These arguments also provide an algebraic
counterpart to the Banach lattice amalgamation methods in \cite{amalg}.
The Archimedean category behaves differently: amalgamation can fail, even
though its free products still preserve embeddings
(Proposition~\ref{prop:injectivity-avl-freeprod}). Thus stability under free
products and stability under general push-outs must be distinguished.

Our first main result characterizes the finitely generated Archimedean
vector lattices $X$ for which $X*_{\scriptscriptstyle{\mathcal{VL}}}Y$ is Archimedean for every finitely generated Archimedean $Y$. We prove that this property is equivalent to projectivity of $X$
(Proposition~\ref{prop:X-fg-projective}). Equivalently, $X$ is finitely
presented, embeds into a free vector lattice, or is represented by a
polyhedral cone in the Baker--Beynon duality (see \cite{baker,Beynon1975}). The latter equivalences are
classical, but the free-product condition gives a further characterization of
the same class. The proof uses the description of finitely generated
Archimedean vector lattices as lattices of piecewise linear functions on
closed cones, under which Archimedean free products correspond to products
of cones. For every non-polyhedral cone, one can produce a partner for which the
ordinary free product is non-Archimedean. For arbitrary Archimedean $X$,
the free product with every Archimedean $Y$ is Archimedean precisely when
every finitely generated sublattice of $X$ is projective
(Corollary~\ref{cor:locally-projective-freeprod}).
We also show that this local condition does not imply that $X$ itself
embeds into a free vector lattice. These results provide a connection between preservation of the Archimedean property, finite presentation and polyhedral geometry.

For Banach lattice factors $(X_i)_{i\in I}$, we identify precisely the
algebraic sublattice generated by their canonical copies in the Banach
lattice free product. The canonical maps factor as
\[
 \underset{i\in I}{\bigfree\nolimits_{\scriptscriptstyle{\mathcal{VL}}}}X_i
 \twoheadrightarrow
 \underset{i\in I}{\bigfree\nolimits_{\scriptscriptstyle{\mathcal{AVL}}}}X_i
 \hookrightarrow
 \underset{i\in I}{\bigfree\nolimits_{\scriptscriptstyle{\mathcal{BL}}}}X_i,
\]
where the first map is Archimedeanization and the second is injective
(Theorem~\ref{theo:avl-into-bl}). Moreover, the Banach lattice free product
is the completion of the Archimedean free product under the largest lattice
norm making the canonical maps contractive. Thus Archimedeanization
accounts for all the algebraic identifications required before taking
the norm completion. This step is necessary: the
vector-lattice free product of two free Banach lattices can fail to be
Archimedean (see Example~\ref{ex:non-archimedean-freeprod}).

The join representation is also useful beyond factors of the form $C(K)$.
We first prove that injective lattice homomorphisms from $C(K)$ and $C(L)$
into arbitrary Banach lattices induce an injective map on their free
products (Theorem~\ref{theo:ck-injectivity}). We then show that every element
of the Banach lattice free product of two Banach lattices comes from the free product of
two principal ideals, equipped with their order-unit norms
(Proposition~\ref{prop:any-element-range-principal-ideals}). This brings the
representation by continuous functions into the study of individual elements of a general free
product. It implies that Archimedean free products of order dense
sublattices are order dense in the corresponding two-factor Banach lattice free product
(Corollary~\ref{cor:order-denseness}), and that each factor is a regular
sublattice of an arbitrary Banach lattice free product
(Corollary~\ref{cor:free-factor-regular}). We give a countable family for
which the algebraic sublattice generated by the factors is not order
dense, despite being norm dense. Finally, every uniformly bounded family
of injective Banach lattice homomorphisms induces an injective
homomorphism between the free products
(Proposition~\ref{prop:injectivity-bl-freeprod}). This extends the embedding results of
\cite{freeprods_BL} and clarifies how injectivity and order approximation
behave under completion.

Our final group of results concerns extension problems for vector
lattices. Motivated by the connection between amalgamation, relative
injectivity and universal disposition in separable Banach lattices developed in \cite{amalg},
we consider the class
\[
 \mathcal C=\{X:|\operatorname{Spec}_p(X)|\leq\aleph_0\},
\]
where $\operatorname{Spec}_p(X)$ denotes the set of proper prime ideals of a vector lattice $X$. A vector lattice is called $\mathcal C$-injective, or countable-prime-injective, if lattice
homomorphisms into it extend from sublattices of members of $\mathcal C$.
Universal disposition for $\mathcal C$ requires embeddings to admit
extensions that are also embeddings.

Writing $\mathfrak c=2^{\aleph_0}$, we prove that every member of
$\mathcal C$ has cardinality at most $\mathfrak c$, although it may have uncountable
dimension (Proposition~\ref{cpi:cardinality} and
Example~\ref{cpi:lex-product}). This bound permits a transfinite push-out
construction which embeds every vector lattice into one of universal
disposition for $\mathcal C$ (Theorem~\ref{cpi:ud-existence}). Such targets
are countable-prime-injective, and there exists a nonzero target of
universal disposition of cardinality at most $2^{\mathfrak c}$.
The construction takes place in $\mathcal{VL}$, since $\mathcal C$ is
closed under sublattices and quotients but not under free products or
countable increasing unions.

A different construction gives smaller vector lattices which are $\mathcal C$-injective. For every
nonprincipal ultrafilter $\mathcal U$ on $\mathbb N$, the ordered
ultrapower $\mathbb R^{\mathbb N}/\mathcal U$ is countable-prime-injective
and has cardinality $\mathfrak c$ (Theorem~\ref{cpi:ultrapower}). More
generally, countable interpolation suffices for a totally ordered real
vector space to have this extension property. The proof relates prime
quotients and Archimedean classes to the realization of countable order
cuts in ultrapowers. Since a totally ordered vector lattice cannot
contain coordinatewise ordered $\mathbb R^2$, these examples also separate
relative injectivity from universal disposition. On the other hand, every
nonzero $\mathcal C$-injective vector lattice has at least $\mathfrak c$ proper
prime ideals and is necessarily non-Archimedean. The obstruction persists
even when all vector lattices involved are Archimedean: no nonzero Archimedean vector
lattice is injective relative to $\mathcal C\cap\mathcal{AVL}$
(Theorem~\ref{cpi:archimedean-obstruction}). This exhibits a limitation of
purely order-theoretic extension properties which is relevant to their
comparison with extension problems in Banach lattices.

The paper is organized as follows. Section~\ref{sec:vl-preliminaries}
collects the algebraic constructions and the prime-ideal arguments.
In Section~\ref{sec:ck-free-products}, we develop the $C(K)$ injectivity result.
Section~\ref{sec:archimedean-free-products} studies Archimedeanization and
the geometric criterion for free products, and
Section~\ref{sec:banach-free-products} establishes the completion and
order-theoretic consequences for Banach lattices. Finally, Section~\ref{sec:countable-spectrum} deals with universal
disposition and relative injectivity with respect to the class of vector lattices with countable prime spectra.

\subsection{Formalization}
The results and numbered definitions in this paper have been formalized in
Lean~4.35.0. The development uses
\href{https://github.com/leanprover-community/mathlib4/tree/56a005855ea17e0b0958f73ec801ae140b0a1679}{Mathlib~4.35.0}
and the Banach lattice library
\href{https://github.com/davidmunozlahoz/banlat/tree/c3c5e1bfc4b2ae61d911aa63d331b9b05a5107d6}{BanLat},
pinned to commits \texttt{56a005855ea1} and \texttt{c3c5e1bfc4b2}, respectively.
The formalization is maintained in the
\href{https://github.com/GonzMA/FreeProductsVecLat/tree/\leanrevision}{FreeProductsVecLat}
repository; the version linked here is revision \texttt{ca3135e88ade}.

Each numbered statement carries a linked $\boldsymbol{\forall}$ symbol and
the name of its corresponding Lean declaration. When a statement is split
into several formal results, we provide the relevant links together.
The links target source lines at the fixed revision above, independently
of subsequent changes to the repository. The mathematical sources contain
no proof placeholders, and the build includes a transitive check that the
registered results use only Lean's standard axioms.

\subsection{Statement on AI use}

GPT models were used in developing several of the mathematical arguments and
constructions in this paper. GPT-5.5 was used to develop the proof, based on
prime ideals, that push-outs of vector lattices preserve injective lattice
homomorphisms (Proposition~\ref{prop:pushout-vl}). The same model was used
for the proof that injective lattice homomorphisms from $C(K)$ spaces induce
an injective map on their Banach lattice free products
(Theorem~\ref{theo:ck-injectivity}), and for the proof that every element of
the Banach lattice free product of two Banach lattices belongs to the image
of the free product of two principal ideals equipped with their order-unit
norms (Proposition~\ref{prop:any-element-range-principal-ideals}).

GPT-5.6 was used to develop the discussion of polyhedral cones and the main
characterization result in Section~\ref{sec:archimedean-free-products}
(Proposition~\ref{prop:X-fg-projective} and
Corollary~\ref{cor:locally-projective-freeprod}). In particular, it suggested
using cone duality and helped formulate the precise conditions under which
vector lattice free products of Archimedean factors remain Archimedean.

GPT-6 Astra was used to develop some of the results and proofs in Section~\ref{sec:countable-spectrum}, on vector lattices with countable
spectra, relative injectivity, and universal disposition.

GPT-6 Astra and Sol were also used to assist with the Lean formalization.

The authors subsequently reviewed, corrected, and refined the AI-assisted
arguments, adapting them to the mathematical framework and exposition of
the paper.

\section{Preliminaries on vector lattices}\label{sec:vl-preliminaries}

Throughout, all vector spaces and lattices are over the real field. An
ideal of a vector lattice is a vector subspace which issolid ($x\in I$ whenever $|x|\leq |y|$ for some $y\in I$). Lattice homomorphisms are linear and
preserve finite suprema and infima. We write
$\lat(A)$ for the vector sublattice generated by a set $A$.
For background on vector and Banach lattices, see \cite{luxemburg,nieberg}.

For a vector space $E$, the free vector lattice $\operatorname{FVL}[E]$
is generated by a linear copy of $E$ and has the property that every
linear map from $E$ into a vector lattice extends uniquely to a lattice
homomorphism; see \cite{baker}.

When $E$ is a Banach space, $\operatorname{FBL}[E]$ denotes the free
Banach lattice generated by $E$. It contains a canonical linear
isometric copy of $E$, and every bounded linear map from $E$ into a
Banach lattice extends uniquely to a lattice homomorphism from
$\operatorname{FBL}[E]$ of the same norm. This construction was
introduced by Avil\'es, Rodr\'iguez and Tradacete~\cite{freebasico},
generalizing the construction over sets due to de Pagter and
Wickstead~\cite{projfree}. For a systematic study of their structure
and their relation to the geometry of the generating Banach space,
see Oikhberg, Taylor, Tradacete and Troitsky~\cite{freebanlat}.

For free products of Banach lattices, we follow the conventions
of~\cite{freeprods_BL}.

Given a family of vector lattices $(X_i)_{i\in I}$, a free product of this
family is a vector lattice $X$ together with lattice homomorphisms
$\varphi_i:X_i\longrightarrow X$ such that, for every family of lattice
homomorphisms $T_i:X_i\longrightarrow Y$, there is a unique lattice
homomorphism $T:X\longrightarrow Y$ satisfying $T_i=T\varphi_i$ for every
$i\in I$. The corresponding definitions in $\mathcal{AVL}$ and
$\mathcal{BL}$ use Archimedean targets and contractive Banach lattice
homomorphisms, respectively. By rescaling, a uniformly bounded family
of Banach lattice homomorphisms also induces a bounded homomorphism
on the Banach lattice free product, of norm at most $\sup_i\|T_i\|$.

\begin{prop}
\leanref{FreeProducts/VectorLattice/Basic.lean}{121}{nonempty_freeProduct}
The free product of a family of vector lattices always exists.
\end{prop}
\begin{proof}
    Let $(X_i)_{i\in I}$ be a family of vector lattices, and let $\iota_i:X_i\longrightarrow Z=(\bigoplus_{j\in I} X_j)$ be the inclusions into their algebraic direct sum. Given lattice homomorphisms $T_i:X_i\longrightarrow Y$, set $S=\oplus_i T_i:Z\longrightarrow Y$. The induced lattice homomorphism $\hat{S}:\operatorname{FVL}[Z]\longrightarrow Y$ satisfies $\hat{S}\delta=S$, where $\delta:Z\longrightarrow \operatorname{FVL}[Z]$ is the canonical inclusion.

    Consider the ideal $J$ of $\operatorname{FVL}[Z]$ generated by elements of the form $\delta\iota_i|x|-|\delta\iota_i(x)|$, with $i\in I$ and $x\in X_i$. If $q:\operatorname{FVL}[Z]\longrightarrow X=\operatorname{FVL}[Z]/J$ denotes the quotient map and $\varphi_i=q\delta\iota_i$, then 
    \[
    \varphi_i|x|=q\delta\iota_i|x|=q|\delta\iota_i(x)|=|\varphi_i(x)|.
    \]

$\hat{S}$ vanishes on $J$ because 
\[
\hat{S}(\delta\iota_i|x|-|\delta\iota_i(x)|)=S\iota_i|x|-|\hat{S}\delta\iota_i(x)|=T_i|x|-|T_i(x)|=0.
\] 
Hence $\hat{S}$ induces a lattice homomorphism $T:X\longrightarrow Y$ with $Tq=\hat{S}$. The map $T$ satisfies 
\[
T\varphi_i=Tq\delta\iota_i=\hat{S}\delta\iota_i=S\iota_i=T_i.
\] 
Uniqueness follows because the sublattices $\varphi_i(X_i)$ generate $X$.
\end{proof}

As usual, the free product of a family of vector lattices is unique up to lattice isomorphism. We denote by $\underset{i\in I}{\bigfree\nolimits_{\scriptscriptstyle{\mathcal{VL}}}} X_i$ the vector lattice free product of the family $(X_i)_{i\in I}$.

\begin{obs}
\leanref{FreeProducts/VectorLattice/Basic.lean}{172}{inclusion_injective}
The inclusions $\varphi_{i_0}:Y_{i_0}\rightarrow \bigfree\nolimits_{\scriptscriptstyle{\mathcal{VL}}} Y_i$ are injective.
\end{obs}
\begin{proof}
Choose $T_{i_0}=\operatorname{id}_{Y_{i_0}}$ and $T_i=0$ when $i\neq i_0$. Then $(\bigfree_{i} T_i)\varphi_{i_0}=\operatorname{id}_{Y_{i_0}}$, so $\varphi_{i_0}$ is injective.
\end{proof}

\begin{prop}
\leanref{FreeProducts/VectorLattice/Basic.lean}{184}{isFreeOver_directSum}
For every family of vector spaces $(X_i)$, \[\operatorname{FVL}[\oplus X_i]\simeq \bigfree\nolimits_{\scriptscriptstyle{\mathcal{VL}}} \operatorname{FVL}[X_i].\]
\end{prop}

\begin{proof}
The argument in \cite[Proposition 5.5]{freeprods_BL} shows that $\operatorname{FVL}[\oplus X_i]$ has the required universal property.
\end{proof}

\begin{prop}\label{prop:quotients-freeprod}
\leanref{FreeProducts/VectorLattice/Basic.lean}{257}{exists_quotient_equiv}
Let $J_i\subset X_i$ be ideals. Then
\[
\bigfree\nolimits_{\scriptscriptstyle{\mathcal{VL}}} (X_i/J_i)\simeq \big(\bigfree\nolimits_{\scriptscriptstyle{\mathcal{VL}}} X_i\big )/J,
\]
where $J$ is the ideal generated by the $\varphi_{X_i}(J_i)$.
\end{prop}

\begin{proof}
Let $q_i:X_i\to X_i/J_i$ and $q:\bigfree\nolimits_{\scriptscriptstyle{\mathcal{VL}}} X_i\to Z:=(\bigfree\nolimits_{\scriptscriptstyle{\mathcal{VL}}} X_i)/J$ be the quotient maps. The vector lattice $(\bigfree\nolimits_{\scriptscriptstyle{\mathcal{VL}}} X_i)/J$, together with the homomorphisms $\psi_i$ given by $\psi_iq_i=q\varphi_i$, is a free product of the family $(X_i/J_i)$:

Given lattice homomorphisms $T_i:X_i/J_i\to Y$, the induced homomorphism vanishes on each $\varphi_i (J_i)$, and hence on $J$. The induced homomorphism
\[
T:=\widetilde{\bigfree_i (T_iq_i)}:Z\to Y
\]
satisfies $T\psi_i=T_i$ and is unique with this property.
\end{proof}

\begin{prop}
\leanref{FreeProducts/VectorLattice/Basic.lean}{381}{nonempty_pushout}
The category of vector lattices with lattice homomorphisms has push-outs. In particular, given lattice homomorphisms
\[
T_1:X_0\to X_1,\qquad T_2:X_0\to X_2,
\]
there exist a vector lattice $PO$ and lattice homomorphisms $S_i$, $i=1,2$, such that the diagram
\[
\xymatrix{X_1 \ar[r]^{S_1} & PO \\
X_0 \ar[u]^{T_1} \ar[r]^{T_2} & X_2 \ar[u]_{S_2}}
\]
commutes and, for every vector lattice $Y$ and homomorphisms $S'_i:X_i\to Y$ making the corresponding square commute, there is a unique lattice homomorphism $R:PO\to Y$ such that $S_i'=RS_i$.
\end{prop}

\begin{proof}
Let $V=X_1*_{\scriptscriptstyle{\mathcal{VL}}}X_2$, with canonical
maps $\varphi_1,\varphi_2$, and let $J$ be the ideal generated by
$\varphi_1T_1(x)-\varphi_2T_2(x)$ for $x\in X_0$. Set $PO=V/J$ and
$S_i=q\varphi_i$, where $q:V\to V/J$ is the quotient map. These maps
satisfy $S_1T_1=S_2T_2$. For any compatible pair $S'_i:X_i\to Y$, the
free-product homomorphism $V\to Y$ vanishes on $J$, so it factors
uniquely through $PO$. This proves the universal property.
\end{proof}

\subsection{Totally ordered vector spaces and prime ideals}
\begin{defin}
\leanref{FreeProducts/VectorLattice/PrimeIdeals.lean}{23}{IsPrime}
A prime ideal $P$ of a vector lattice $X$ is an ideal such
that, whenever $a,b\in X_+$ satisfy $a\wedge b=0$, either $a\in P$ or
$b\in P$. 
\end{defin}

In this section we allow $P=X$. The
prime spectra in Section~\ref{sec:countable-spectrum} contain only
proper prime ideals.

\begin{prop}
\leanref{FreeProducts/VectorLattice/PrimeIdeals.lean}{27}{isTotal_quotient_iff_isPrime}
If $J$ is an ideal of $X$, then $X/J$ is totally ordered if and only if $J$ is prime.
\end{prop}
\begin{proof} See \cite[Theorem~3.1]{JohnsonKist1962}.
\end{proof}

\begin{prop}\label{prop:prime-ideals-extcon}
\leanref{FreeProducts/VectorLattice/PrimeIdeals.lean}{49}{exists_prime_restriction}
Suppose $Y\subset X$ is a vector sublattice and $P\subset X$ is a prime ideal. Then $P\cap Y$ is a prime ideal of $Y$. \leanref{FreeProducts/VectorLattice/PrimeIdeals.lean}{155}{exists_prime_extension}
Conversely, if $P\subset Y$ is a prime ideal of $Y$, there is a prime ideal of $X$, denoted by $Q$, such that $Q\cap Y=P$.
\end{prop}
\begin{proof}
Restriction is immediate. For extension, let $I$ be the ideal generated
by $P$ in $X$. Explicitly,
$I=\{x\in X:|x|\leq p\text{ for some }p\in P_+\}$, so $I\cap Y=P$.
Passing to $X/I$ and $Y/P$ reduces
the problem to the case where $Y$ is totally ordered and $P=\{0\}$.
By Zorn's lemma, choose an ideal $Q$ maximal subject to $Q\cap Y=\{0\}$.
If $Q$ were not prime, there would be disjoint $a,b\in X_+\setminus Q$.
Maximality gives nonzero $y_a,y_b\in Y_+$, an element $q\in Q_+$ and
$\lambda>0$ such that
\[
y_a\leq q+\lambda a,\qquad y_b\leq q+\lambda b.
\]
Since $Y$ is totally ordered,
\[
0<y_a\wedge y_b\leq(q+\lambda a)\wedge(q+\lambda b)=q,
\]
contradicting $Q\cap Y=\{0\}$. Thus $Q$ is prime.
\end{proof}

The following amalgamation result for totally ordered vector spaces is a special case of the more general results of Pierce~\cite{Pierce:1973,Pierce:1976} and Hager and Madden~\cite{HagerMadden1983}. We include a direct proof, which is also used in the proof of Proposition~\ref{prop:pushout-vl}.

\begin{prop}\label{prop:TOVS-square}
\leanref{FreeProducts/VectorLattice/TotalAmalgamation.lean}{195}{nonempty_totalAmalgam}
Let \(X\), \(Y\), \(Z\) be totally ordered vector spaces, and suppose there are injective lattice homomorphisms
\(\iota_X:Z\to X\) and \(\iota_Y:Z\to Y\), which identify $Z$ with ordered subspaces of $X$ and $Y$. Then there exist a totally ordered vector space
\(Z_0\) and injective lattice homomorphisms \(\alpha:X\to Z_0\) and
\(\beta:Y\to Z_0\) such that \[
\alpha\iota_X=\beta\iota_Y.
\]
\end{prop}

\begin{proof}
Consider the vector-space amalgam
\[
Z_0=(X\oplus Y)/\{(\iota_X(c),-\iota_Y(c)):c\in Z\},
\]
with canonical injective maps $\alpha$ and $\beta$, which agree on $Z$.
The wedge $C_0=\alpha(X_+)+\beta(Y_+)$ is a cone. Indeed, if
$\alpha(a)+\beta(b)=0$ with $a,b\geq0$, then
$a=\iota_X(c)$ and $b=-\iota_Y(c)$ for some $c\in Z$; hence $c\geq0$
and $c\leq0$, so $a=b=0$. Now let $w\in C_0\cap(-C_0)$ and write
$w=\alpha(a_1)+\beta(b_1)$ and $-w=\alpha(a_2)+\beta(b_2)$, where
$a_1,a_2\in X_+$ and $b_1,b_2\in Y_+$. Adding gives
\[
\alpha(a_1+a_2)+\beta(b_1+b_2)=0.
\]
The preceding argument yields $a_1+a_2=0$ and $b_1+b_2=0$. Positivity then
forces all four elements to be zero, so $w=0$. Thus
$C_0\cap(-C_0)=\{0\}$.

By Zorn's lemma, extend $C_0$ to a maximal cone $M$ in $Z_0$.
It is total: if $v\notin M\cup(-M)$, then $M+\mathbb R_+v$ is a
strictly larger cone, since
$m+n+(\lambda+\mu)v=0$, with $m,n\in M$ and $\lambda,\mu\geq0$,
forces $\lambda=\mu=0$ and $m=n=0$.
The order defined by $M$ therefore makes $Z_0$ a totally ordered vector space.
The maps $\alpha$ and $\beta$ preserve positivity and, being injective
maps between totally ordered spaces, also reflect it and preserve lattice
operations.
\end{proof}
\subsection{Embeddings and colimits of vector lattices}
The following classical fact is discussed by Hager and Madden
\cite[pp.~183--184]{HagerMadden1983}, who also note that their argument
for abelian lattice-ordered groups applies to vector lattices.
We include a direct proof using prime ideals, following the method of
\cite[Lemma~2 and the proof of Theorem~1(a)]{Pierce:1976} for abelian
lattice-ordered groups. This presentation also makes explicit the parallel
with the proof of \cite[Theorem~4.4]{amalg}: prime quotients and totally ordered
vector lattices play the role of the Hahn--Banach theorem and the $L_1$-representations used there.

\begin{prop}\label{prop:pushout-vl}
\leanref{FreeProducts/VectorLattice/Amalgamation.lean}{92}{left_injective}
Let \(X_0\), \(X_1\), and \(X_2\) be vector lattices, and let \[
\xymatrix{X_1 \ar[r]^{S_1} & PO \\
X_0 \ar[u]^{T_1} \ar[r]^{T_2} & X_2 \ar[u]^{S_2}}
\]
be a push-out square. If $T_2$ is injective, then so is $S_1$.
\end{prop}
\begin{proof}
Identify $X_0$ with its image under $T_2$ in $X_2$. Fix $0\neq x\in X_1$
and choose a prime ideal $P_1$ with $x\notin P_1$.
The inverse image $P_0=T_1^{-1}(P_1)$ is prime in $X_0$.
By Proposition \ref{prop:prime-ideals-extcon}, choose a prime ideal
$P_2$ in $X_2$ such that $P_2\cap X_0=P_0$.
The induced maps
\[
\widetilde T_1:X_0/P_0\to X_1/P_1,\qquad
\widetilde T_2:X_0/P_0\to X_2/P_2
\]
are injective lattice homomorphisms between totally ordered vector spaces.
By Proposition \ref{prop:TOVS-square}, they admit an amalgamation
into a totally ordered vector space $Z$. Writing $q_0,q_1,q_2$ for the
quotient maps, we obtain the commutative diagram
\[
\xymatrix{& X_1 \ar[r]^-{q_1} & X_1/P_1 \ar[rd]^{R_1} \\
X_0 \ar[ru]^{T_1} \ar[rd]^{T_2} \ar[r]^{q_0} & X_0/P_0 \ar[ru]^{\widetilde T_1} \ar[rd]^{\widetilde T_2} && Z\\
& X_2 \ar[r]^-{q_2} & X_2/P_2 \ar[ru]^{R_2}}
\]
with $R_1$ and $R_2$ injective.
The push-out property gives $R:PO\to Z$ with $RS_1=R_1q_1$.
Consequently,
\[
RS_1(x)=R_1q_1(x)\neq0,
\]
so $S_1$ is injective.
\end{proof}

\begin{cor}\label{cor:freeprod-vl-injective}
\leanref{FreeProducts/VectorLattice/Amalgamation.lean}{207}{induced_injective}
Let $T_i:Y_i\to X_i$ be injective vector lattice homomorphisms, with $i\in I$. The free product map
\[\overline{\bigfree} T_i:\bigfree\nolimits_{\scriptscriptstyle{\mathcal{VL}}} Y_i\to \bigfree\nolimits_{\scriptscriptstyle{\mathcal{VL}}} X_i\]
is injective.
\end{cor}

\begin{proof}
Identify each $Y_i$ with its image in $X_i$. For two factors, replacing
one factor by the larger lattice gives a push-out of its inclusion
along its canonical map into the free product. Proposition~\ref{prop:pushout-vl}
makes the resulting map injective. Replacing the factors successively
proves the finite case. For an arbitrary family, every element of a free
product involves only finitely many factors. The canonical map from a
finite subproduct into the full product has a left inverse obtained by
sending all remaining factors to zero, so it is injective. The finite
case therefore proves the assertion for arbitrary $I$.
\end{proof}

\section{Revisiting free products of \texorpdfstring{$C(K)$}{C(K)}}\label{sec:ck-free-products}

In this section we prove that injective lattice homomorphisms from $C(K_1)$ and $C(K_2)$ into arbitrary Banach lattices induce an injective map on their free products. The proof uses the join representation of free products of $C(K)$ spaces given in \cite{freeprods_BL}, together with a simple approximation argument.

\begin{lema}\label{lema:X_U}
Let \(T:C(K)\to X\) be an injective lattice homomorphism into a Banach lattice $X$. Let \(k_0\in K\), and let \(U\) be
an open neighborhood of \(k_0\) in $K$. Let \(I_U\) be the closed ideal of \(X\)
generated by $T\{f:f|_U=0\}$, and let
\(q:X\to X_U=X/I_U\) be the quotient lattice homomorphism. \leanref{FreeProducts/Consequences/ContinuousFunctionEmbeddings.lean}{91}{local_zero_ideal_unit_not_mem}
Then
\(e=qT(1_K)\) is nonzero. \leanref{FreeProducts/Consequences/ContinuousFunctionEmbeddings.lean}{112}{local_quotient_estimate}
Moreover, if \(f\in C(K)\) and
\(\varepsilon>0\) satisfy \(|f-f(k_0)1_K|\leq\varepsilon 1_K\) on \(U\), then \[
\bigl|qT(f)-f(k_0)e\bigr|\leq \varepsilon e .
\]
\end{lema}

\begin{proof}
Choose \(0<g\in C(K)\) such that
\(\operatorname{supp} g\subset U\).
If \(f\in C(K)\) and \(f|_U=0\), then \(|f|\wedge g=0\). Since \(T\) is
a lattice homomorphism, \(|T(f)|\wedge T(g)=0\). Hence every generator
\(T(f)\) of \(I_U\) belongs to the disjoint complement of \(T(g)\), and
the disjoint complement of \(T(g)\) is a closed ideal. Therefore \(I_U\)
is contained in the disjoint complement of \(T(g)\).

If \(qT(g)=0\), then \(T(g)\in I_U\), so \(T(g)\perp T(
g)\) and hence
\(T(g)=0\). This contradicts the injectivity of \(T\). Thus
\(qT(g)\neq 0\). Since \(0\leq T(g)\leq \|g\|_\infty T(1_K)\), we have that 
\[ 
0\neq qT(g)\leq \|g\|_\infty qT(1_K)=\|g\|_\infty e,
\]
so \(e\neq 0\).
\end{proof}

For convenience, we now recall the representation of free products of $C(K)$ spaces
from \cite[Theorem 6.2]{freeprods_BL}. First, let us recall the definition of the join. For nonempty compact Hausdorff spaces $K_1$ and $K_2$, their join is
\[
K_1*K_2=(K_1\times K_2\times[0,1])/\sim,
\]
where $(k,l,0)\sim(k,l',0)$ and $(k,l,1)\sim(k',l,1)$ for all
$k,k'\in K_1$ and $l,l'\in K_2$. Write $[k,l,t]$ for the equivalence class
of $(k,l,t)$.

\begin{theo}\label{theo:ck-join-representation}
\leanref{FreeProducts/BanachLattice/Join.lean}{821}{exists_join_equiv}
Let $K_1$ and $K_2$ be nonempty compact Hausdorff spaces. There is a lattice isomorphism
\[
C(K_1)*_{\scriptscriptstyle{\mathcal{BL}}}C(K_2)\simeq C(K_1*K_2)
\]
under which the canonical embeddings satisfy
\[
\varphi_1(f)([k,l,t])=(1-t)f(k),\qquad
\varphi_2(g)([k,l,t])=t g(l).
\]
\end{theo}

Note that here $C(K_1*K_2)$ is always equipped with its usual supremum norm. The displayed isomorphism is an isomorphism of Banach lattices, in general it is not an isometry for
the free product norm.

\begin{theo}\label{theo:ck-injectivity}
\leanref{FreeProducts/Consequences/ContinuousFunctionEmbeddings.lean}{443}{induced_injective_of_continuousFunctions}
Let $K_1,K_2$ be nonempty compact Hausdorff spaces and let \(X_1,X_2\) be
Banach lattices. Suppose there are injective lattice
homomorphisms \(T_1:C(K_1)\to X_1\) and \(T_2:C(K_2)\to X_2\). The induced
lattice homomorphism \[T_1* T_2:C(K_1)*_{\scriptscriptstyle{\mathcal{BL}}}C(K_2)\to X_1*_{\scriptscriptstyle{\mathcal{BL}}}X_2
\]
is injective.

\end{theo}
\begin{proof}
Identify the domain with $C(K_1*K_2)$ by the isomorphism from Theorem \ref{theo:ck-join-representation} and write $T=T_1* T_2$.
Fix $p_0=[k_0,l_0,t_0]\in K_1*K_2$. We will show that $T(h)=0$ implies
$h(p_0)=0$ by constructing lattice homomorphisms that approximate evaluation
at $p_0$.

Let $U$ and $V$ be open neighborhoods of $k_0$ and $l_0$. Applying
Lemma \ref{lema:X_U} to $T_1$ and $U$, and to $T_2$ and $V$, gives
quotient lattice homomorphisms $q_1:X_1\to X_{1,U}$ and $q_2:X_2\to X_{2,V}$.
The elements $e=q_1T_1(1_{K_1})$ and $d=q_2T_2(1_{K_2})$ are nonzero and
positive. Put
\[
Z=X_{1,U}\widehat{\otimes}_{|\pi|}X_{2,V},
\]
the Fremlin projective tensor product \cite{Fremlin1974}. Its canonical
bilinear map is a lattice bimorphism and its norm is a cross norm:
$|x\otimes y|=|x|\otimes|y|$ and
$\|x\otimes y\|=\|x\|\|y\|$.
Thus $w=e\otimes d>0$, and the following maps are bounded lattice
homomorphisms $S_1:X_1\to Z$ and $S_2:X_2\to Z$:
\[
S_1(x)=(1-t_0)q_1(x)\otimes d,\qquad
S_2(y)=t_0 e\otimes q_2(y).
\]
The universal property gives $S=S_1*S_2:X_1*_{\scriptscriptstyle{\mathcal{BL}}}X_2\to Z$.
Thus $\Theta=ST:C(K_1*K_2)\to Z$ satisfies
\[
\Theta(\varphi_1(f))=(1-t_0)q_1T_1(f)\otimes d,\qquad
\Theta(\varphi_2(g))=t_0 e\otimes q_2T_2(g).
\]
In particular, since $1_{K_1*K_2}=\varphi_1(1_{K_1})+\varphi_2(1_{K_2})$, we have
\[
\Theta(1_{K_1*K_2})=w.
\]
Compare this map with
$\Theta_0:C(K_1*K_2)\to Z$, defined by $\Theta_0(h)=h(p_0)w$: on the
canonical copies it satisfies
\[
\Theta_0(\varphi_1(f))=(1-t_0)f(k_0)w,\qquad
\Theta_0(\varphi_2(g))=t_0g(l_0)w.
\]
Note these constructions depend on $U$ and $V$, but for the sake of brevity we suppress this dependence in the notation.

We claim that, for each $h\in C(K_1*K_2)$ and $\eta>0$, the neighborhoods
$U$ and $V$ can be chosen so that
\[
|\Theta(h)-\Theta_0(h)|=|\Theta(h)-h(p_0)w|\leq\eta w.
\]
Indeed, first suppose that $h$ belongs to the sublattice generated by the canonical
copies of $C(K_1)$ and $C(K_2)$. It is a finite lattice expression in functions
$f_1,\ldots,f_m\in C(K_1)$ and $g_1,\ldots,g_n\in C(K_2)$.
Given $\delta>0$, choose $U$ and $V$ so that
\[
|f_i-f_i(k_0)1_{K_1}|\leq\delta\quad\text{on }U,\qquad
|g_j-g_j(l_0)1_{K_2}|\leq\delta\quad\text{on }V
\]
for all $i,j$. Lemma \ref{lema:X_U} then gives
\[
|q_1T_1(f_i)-f_i(k_0)e|\leq\delta e,\qquad
|q_2T_2(g_j)-g_j(l_0)d|\leq\delta d.
\]
Tensoring with $d$ or $e$ and multiplying by $1-t_0$ or $t_0$ yields
\[
|\Theta(\varphi_1(f_i))-\Theta_0(\varphi_1(f_i))|\leq\delta w,\qquad
|\Theta(\varphi_2(g_j))-\Theta_0(\varphi_2(g_j))|\leq\delta w.
\]
Using the inequalities
\[
|(u\vee v)-(u'\vee v')|\leq|u-u'|+|v-v'|,\qquad
|(u\wedge v)-(u'\wedge v')|\leq|u-u'|+|v-v'|,
\]
a standard induction argument gives $|\Theta(h)-\Theta_0(h)|\leq C_h\delta w$,
where $C_h$ depends only on the chosen expression, not on $\delta$.
Taking $\delta$ sufficiently small proves the claim on the generated sublattice.

For general $h$, choose $\tilde h\in \lat(\varphi_1(C(K_1)),\varphi_2(C(K_2)))$
with $\|h-\tilde h\|_\infty<\eta/3$. Since both $\Theta$ and $\Theta_0$
send $1_{K_1*K_2}$ to $w$, each maps $h-\tilde h$ to an element of absolute
value at most $(\eta/3)w$. Choose $U,V$ using the preceding case for
$\tilde h$ and $\eta/3$. The triangle inequality now gives
$|\Theta(h)-\Theta_0(h)|\leq\eta w$, proving the claim.

Finally, suppose $T(h)=0$. Then $\Theta(h)=0$ for every choice of $U,V$.
If $h(p_0)\neq0$, choose $0<\eta<|h(p_0)|$ and apply the claim to obtain
\[
|h(p_0)|w=|\Theta(h)-\Theta_0(h)|\leq\eta w,
\]
contradicting $w>0$. Thus $h(p_0)=0$; since $p_0$ was arbitrary,
$h=0$ and $T$ is injective.
\end{proof}

\section{Archimedean and non-Archimedean free products}\label{sec:archimedean-free-products}
Recall that a vector lattice is Archimedean if $0\leq nx\leq y$ for
every positive integer $n$ implies $x=0$. An ideal is uniformly closed
if it contains the limits of its relatively uniformly convergent
sequences. Here $x_n\to x$ relatively uniformly means that
$|x_n-x|\leq\varepsilon_n u$ for some $u\geq0$ and some real sequence
$\varepsilon_n\downarrow0$. For an ideal $J$ of a vector lattice $X$,
the quotient $X/J$ is Archimedean if and only if $J$ is uniformly
closed; see \cite[Theorem~60.2]{luxemburg}.

Let $\mathfrak a(X)$ be the intersection of all uniformly closed ideals
of $X$. This ideal is uniformly closed, so $X/\mathfrak a(X)$ is
Archimedean. If $Y$ is Archimedean and $T:X\to Y$ is a lattice
homomorphism, then $\ker T$ is uniformly closed and therefore contains
$\mathfrak a(X)$. Thus $T$ factors uniquely as
\[
\xymatrix{X/\mathfrak a(X) \ar@{-->}[rd]^{\tilde T} \\ X \ar[r]^T \ar[u]^q & Y.}
\]
Consequently, $\mathcal{AVL}$ is a reflective subcategory of
$\mathcal{VL}$ \cite[Chapter~IV, \S3]{MacLane1998}. We call $X/\mathfrak a(X)$ the Archimedeanization of
$X$ and denote it by $X_{\mathfrak a}$.

The subcategory of Archimedean vector lattices also has free products and push-outs. These are precisely the Archimedeanizations of the corresponding colimits in $\mathcal{VL}$. As the following example shows, colimits in $\mathcal{VL}$ need not agree with the corresponding colimits in $\mathcal{AVL}$. This contrasts with the fact that free objects in $\mathcal{VL}$ are Archimedean.

The following example is essentially due to K. Pierce \cite{Pierce:1973,Pierce:1976} in the setting of $\ell$-groups.

\begin{ex}
\leanref{FreeProducts/Counterexamples/NonArchimedeanPushout.lean}{45}{archimedean_amalgam_kills_extra_coordinate}
The category $\mathcal{AVL}$ of Archimedean vector lattices and lattice homomorphisms does not have the amalgamation property: two injective lattice homomorphisms with a common domain need not admit a commutative square whose remaining maps are both injective. \leanref{FreeProducts/Counterexamples/NonArchimedeanPushout.lean}{59}{exists_nonarchimedean_pushout}
In particular, the canonical maps into their push-out need not be injective. Take, for example, $X_0=\{(x_n)\in\mathbb{R}^\mathbb{N}:(x_n)\text{ is eventually constant}\}$, $X_1=\mathbb{R}^\mathbb{N}$ and $X_2=X_1\oplus \mathbb{R}$. Let $j_1:X_0\rightarrow X_1$ be the inclusion, and define $j_2:X_0\rightarrow X_2$ by $j_2((x_n))=((x_n),\lim x_n)$. For any Archimedean vector lattice $Z$ and lattice homomorphisms $T_i:X_i\rightarrow Z$ making the diagram
\[
\xymatrix{ X_2 \ar[r]^{T_2} & Z \\ X_0 \ar[r]^{j_1} \ar[u]^{j_2} & X_1 \ar[u]_{T_1}}
\]
commute, we have $T_2(0,1)=0$. Indeed, for each $n\in\mathbb{N}$,
\begin{align*}
    T_1(1,2,3,\ldots)&\geq T_1(1,2,\ldots,n,n,\ldots)=T_1j_1(1,2,\ldots,n,n,\ldots)=T_2j_2(1,2,\ldots,n,n,\ldots)\\
    &=T_2((1,2,\ldots,n,n,\ldots),n)\geq T_2(0,n)=nT_2(0,1).
\end{align*}

Since $Z$ is Archimedean, it follows that $T_2(0,1)=0$. In particular, this holds when $Z$ is the Archimedean push-out of $j_1$ and $j_2$.
\end{ex}

The preceding example shows that the $\mathcal{VL}$ push-out of
Archimedean vector lattices need not be Archimedean. We now study the
corresponding question for free products.

In what follows, we use Baker's conventions from \cite{baker}. A closed cone in $\mathbb R^n$ contains the
origin and is invariant under multiplication by nonnegative scalars, though
it need not be convex. A closed polyhedral cone is a finite union of
closed convex polyhedral cones. An open polyhedral cone is a finite
union of nonempty finite intersections of open half-spaces whose
boundaries pass through the origin. For a function $f$, its open support is
$S(f)=\{x:f(x)\neq0\}$. In the standard representation of free vector lattices,
$\operatorname{FVL}[\mathbb R^n]$ consists of continuous positively
homogeneous piecewise linear functions on $\mathbb R^n$, and their
open supports are precisely the open polyhedral cones.

For finitely generated Archimedean vector lattices, the $\mathcal{AVL}$ free product admits a geometric representation. For a subset $A$ of a vector lattice, let $I(A)$ denote the ideal generated by $A$.
For a closed cone $C\subseteq\mathbb R^n$, write
\[
\Van(C)=\{f\in \operatorname{FVL}[\mathbb{R}^{n}]:f|_C=0\}
\]
for the ideal of functions vanishing on $C$.

For a closed cone \(C\subseteq\mathbb R^m\), put

\[
  \mathscr{L}(C)
  =
  \{f|_C:f\in\operatorname{FVL}[\mathbb{R}^{m}]\}.
\]

 If \(C\subseteq\mathbb R^m\) and
\(D\subseteq\mathbb R^n\) are closed cones, an \emph{\(\ell\)-map}
\(q:C\to D\) is a map \(q=(q_1,\ldots,q_n)\) whose coordinate functions
\(q_i\) belong to \(\mathscr{L}(C)\). In particular, its coordinates
extend to continuous positively homogeneous piecewise linear functions
on the whole ambient space $\mathbb R^m$.
Let \(\mathbf{Cone}^{\ell}\) denote the category of
closed cones in finite-dimensional Euclidean spaces and \(\ell\)-maps.

\begin{theo}[Baker--Beynon duality]\label{theo:baker-beynon-duality}
\leanref{FreeProducts/AVLvsVL/Classical/Duality.lean}{163}{coneDuality}
Let \(\mathcal{AVL}_{\mathrm{fg}}\) denote the category of finitely
generated Archimedean real vector lattices and vector lattice
homomorphisms. There is a contravariant equivalence
\[
  \bigl(\mathbf{Cone}^{\ell}\bigr)^{\mathrm{op}}
  \simeq \mathcal{AVL}_{\mathrm{fg}},
  \qquad C\longmapsto\mathscr L(C),
  \qquad (q:C\to D)\longmapsto
  \bigl(q^*:f\mapsto f\circ q\bigr).
\]
\end{theo}

\begin{proof}
The representation $X\cong\mathscr L(C)$ is
\cite[Theorem 4.3]{baker}. The explicit correspondence on morphisms is in
\cite[Proposition 2.1 and Corollary 1]{Beynon1975}.
\end{proof}

This defines a notion of duality which gives the following description of free products. Let \(X_1,X_2\in\mathcal{AVL}_{\mathrm{fg}}\), and choose closed cones
\(C_1\subseteq\mathbb{R}^{n}\) and
\(C_2\subseteq\mathbb{R}^{m}\), together with corresponding isomorphisms
\[
  X_1\cong\mathscr{L}(C_1)
  \qquad\text{and}\qquad
  X_2\cong\mathscr{L}(C_2).
\]
The free product of \(X_1\) and \(X_2\) in
\(\mathcal{AVL}_{\mathrm{fg}}\) exists and, up to canonical isomorphism,
satisfies
\[
  X_1 *_{\scriptscriptstyle{\mathcal{AVL}_{\mathrm{fg}}}}X_2
  \cong
  \mathscr{L}(C_1\times C_2).
\]
Indeed, $C_1\times C_2$, with the coordinate projections, is the
categorical product in $\mathbf{Cone}^{\ell}$, and the contravariant equivalence
carries products to free products (coproducts).

For two finitely generated Archimedean vector lattices, the $\mathcal{AVL}_{\mathrm{fg}}$ free product agrees with the $\mathcal{AVL}$ free product because the latter is finitely generated.

Define
\[
J_{C_1\times C_2}
:=I\bigl(j_1(\Van(C_1))\cup
          j_2(\Van(C_2))\bigr)
\subseteq \operatorname{FVL}[\mathbb{R}^{n+m}],
\]
where
\[
j_1(f)=f\circ\pi_1,\qquad j_2(g)=g\circ\pi_2
\]
and $\pi_1:\mathbb R^{n+m}\to\mathbb R^n$ and
$\pi_2:\mathbb R^{n+m}\to\mathbb R^m$ denote the coordinate projections. In particular, $J_{C_1\times C_2}\subset \Van(C_1\times C_2)$.

By Theorem \ref{theo:baker-beynon-duality} and Proposition \ref{prop:quotients-freeprod}, for finitely generated Archimedean vector lattices
\[
X_1\simeq \mathscr{L}(C_1)\quad\text{and}\quad X_2\simeq\mathscr{L}(C_2),
\]
we have an isomorphism
\[
X_1*_{\scriptscriptstyle{\mathcal{VL}}} X_2\simeq X_1*_{\scriptscriptstyle{\mathcal{AVL}}} X_2
\]
induced by the identity maps when Archimedeanization identifies no distinct vectors, that is, when the inclusion

\[
J_{C_1\times C_2}\subset \Van(C_1\times C_2)
\]
of ideals in $\operatorname{FVL}[\mathbb R^{n+m}]$ is an equality. We will show that this need not hold in general.

\begin{lema}\label{lema:polyhedral-cover-equivalence}
\leanref{FreeProducts/AVLvsVL/Coproducts/CoveringCriterion.lean}{161}{vanishingIdeal_le_productRelationsIdeal_iff}
Let \(C\subseteq\mathbb R^n\) and \(D\subseteq\mathbb R^m\) be closed
cones. Then
\[
\Van(C\times D)\subset J_{C\times D}
\] if and only if every open polyhedral cone $W\subset \mathbb R^{n+m}$ satisfying
$W\cap(C\times D)=\varnothing$ admits open polyhedral cones \(U\subseteq\mathbb R^n\) and
\(V\subseteq\mathbb R^m\) such that
$
U\cap C=\varnothing$,
$V\cap D=\varnothing$
and
\[
W\subseteq
(U\times\mathbb R^m)\cup(\mathbb R^n\times V).
\]
\end{lema}

\begin{proof}
Let $h\in \Van(C\times D)_+$. We have $h\in J_{C\times D}$ if and only if there are $f\in \Van(C)_+$, $g\in \Van(D)_+$ such that
\[
h\leq f\pi_1+g\pi_2,
\]
or equivalently, by \cite[Theorem 3.4]{baker}, if and only if there exist $f,g$ as above such that
\[
S(h)\subset S(f\pi_1)\cup S(g\pi_2)=(S(f)\times \mathbb R^m)\cup(\mathbb R^n\times S(g)).
\]
Applying the condition to the open polyhedral cone $W=S(h)$ and \cite[Lemma 3.2]{baker}, which identifies open polyhedral cones with open supports of elements of the corresponding
free vector lattices, provide such $f,g$. This proves one implication.

For the converse, assume that $\Van(C\times D)\subset J_{C\times D}$, and let
$W\subset\mathbb R^{n+m}$ be an open polyhedral cone disjoint from
$C\times D$. By \cite[Lemma 3.2]{baker}, there is
$h\in\operatorname{FVL}[\mathbb R^{n+m}]_+$ with $S(h)=W$.
Thus $h\in \Van(C\times D)\subset J_{C\times D}$. By the definition of
the generated ideal, there are $f\in \Van(C)_+$ and $g\in \Van(D)_+$ such that
\[
h\leq f\pi_1+g\pi_2.
\]
Set $U=S(f)$ and $V=S(g)$. These are open polyhedral cones satisfying
$U\cap C=\varnothing$ and $V\cap D=\varnothing$. Since $f,g\geq0$,
\[
W=S(h)\subset S(f\pi_1+g\pi_2)
=(U\times\mathbb R^m)\cup(\mathbb R^n\times V).
\]
Hence $U,V$ satisfy the required covering condition.
\end{proof}

\begin{lema}\label{lema:polyhedral-cone-gives-archimedeanity}
\leanref{FreeProducts/AVLvsVL/Coproducts/CoveringCriterion.lean}{371}{hasPolyhedralCover_of_polyhedral}
If \(C\subseteq\mathbb R^n\) is a closed polyhedral cone and
\(D\subseteq\mathbb R^m\) is any closed cone, then the covering
condition of Lemma \ref{lema:polyhedral-cover-equivalence} holds for the pair
\((C,D)\).
\end{lema}
\begin{proof}
Let \(W\subseteq\mathbb R^{n+m}\) be open polyhedral and disjoint from
\(C\times D\). Define
\[
U=\mathbb R^n\setminus C
\quad\text{and}\quad V=
\pi_2\bigl(W\cap(C\times\mathbb R^m)\bigr).
\]
The cone $U$ is open polyhedral. The set $V$ is polyhedral by elimination
of the $x$ variables in finite systems of homogeneous linear inequalities.
It is open because
$V=\bigcup_{x\in C}\{y:(x,y)\in W\}$ is a union of open sets.
Thus $V$ is also an open polyhedral cone.
The disjointness of \(W\) from \(C\times D\) gives
\(V\cap D=\varnothing\). Finally, if \((x,y)\in W\), then either
\(x\notin C\), in which case \(x\in U\), or \(x\in C\), in which case
\(y\in V\). Therefore
\[
W\subseteq
(U\times\mathbb R^m)\cup(\mathbb R^n\times V).
\]
\end{proof}

For any non-polyhedral cone $C$, we can construct a closed cone \(D_C\) such that the covering condition fails, and hence the corresponding vector lattice free product is not Archimedean.

\begin{lema}\label{lema:cone-bad-partner}
\leanref{FreeProducts/AVLvsVL/Coproducts/NonPolyhedralPartner.lean}{184}{exists_partner_strict_productRelationsIdeal}
Let \(C\subseteq\mathbb R^n\) be a closed nonpolyhedral cone. There is a
closed cone \(D_C\subseteq\mathbb R^{n+1}\) such that

\[
J_{C\times D_C}\subsetneq \Van(C\times D_C).
\]

\end{lema}
\begin{proof}
Define

\[
\Omega_C
=
\left\{
(y,t)\in\mathbb R^{n+1}:
t>0,\
\inf_{x\in C}\|x-y\|_{\infty}<t
\right\}\quad\text{and}\quad
D_C=\mathbb R^{n+1}\setminus\Omega_C.
\]
Since \(C\) is a closed cone, the set \(\Omega_C\) is an open cone and \(D_C\) is a closed cone.

The cone \(\Omega_C\) is not open polyhedral. Indeed, if it were open
polyhedral, its closure would be a closed polyhedral cone: decompose it
into finitely many open convex polyhedral cones and take the finite
union of their closures. Also, linear sections of polyhedral cones are polyhedral. But
\[
\overline{\Omega_C}\cap\big(\mathbb R^n\times\{0\}\big )=C\times \{0\}
\]
is not polyhedral, a contradiction.

Now let
\[
h(x,y,t)
=
\left(
t-\|x-y\|_\infty
\right)^+.
\]
Its open support is the open convex polyhedral cone
\[
W=
\left\{
(x,y,t):
\|x-y\|_\infty< t
\right\}.
\]

By the definition of \(D_C\), $W\cap(C\times D_C)=\varnothing$. Suppose the cover condition held for the pair $(C,D_C)$. Let \(U,V\) be a
corresponding pair of covering open polyhedral cones. For each
\(z=(y,t)\in\Omega_C\), choose \(x\in C\) satisfying \((x,z)\in W\).
Since \(U\cap C=\varnothing\), the cover forces \(z\in V\). Hence
\[
\Omega_C\subseteq V.
\]
On the other hand, \(V\cap D_C=\varnothing\) gives
\(V\subseteq\Omega_C\). Therefore \(V=\Omega_C\), contradicting that
\(\Omega_C\) is not open polyhedral. By Lemma \ref{lema:polyhedral-cover-equivalence},
\[
J_{C\times D_C}\subsetneq \Van(C\times D_C).
\]
\end{proof}

A vector lattice is projective if every lattice homomorphism from it
into a quotient vector lattice lifts through the quotient map. All
references to projectivity below concern the category $\mathcal{VL}$.

\begin{prop} \label{prop:X-fg-projective}
\leanref{FreeProducts/AVLvsVL/Coproducts/Equivalences.lean}{65}{finitely_generated_archimedean_tfae}
Let \(X\) be a finitely generated Archimedean vector lattice. The
following conditions are equivalent:
\begin{enumerate}
\def\labelenumi{\arabic{enumi}.}

\item
  For every finitely generated Archimedean vector lattice \(Y\), the
  ordinary vector-lattice free product \(X*_{\scriptscriptstyle{\mathcal{VL}}}Y\) is
  Archimedean;
\item
  \(X\) is isomorphic to a vector sublattice of a free vector lattice;
\item
  \(X\) is projective;
\item
  \(X\) is finitely presented, that is,
  $X\simeq\operatorname{FVL}[\mathbb R^n]/I$ for some $n$ and some
  finitely generated ideal $I$;
\item
  In any representation \(X\simeq \mathscr L(C)\), the closed cone \(C\)
  is polyhedral.
\end{enumerate}
\end{prop}

\begin{proof}
The equivalence between $(1)$ and $(5)$ follows from Lemma
\ref{lema:cone-bad-partner} and Lemma \ref{lema:polyhedral-cone-gives-archimedeanity}. Also, the equivalence between $(3)$, $(4)$, and $(5)$ is a consequence of \cite[Theorem 5.1 and its proof]{baker}.

For the implication $(3)\implies (2)$, projectivity splits the canonical quotient map
from a finitely generated free vector lattice onto $X$.
For $(2)\implies (5)$, finite generation of $X$ ensures that the image
of an embedding into a free vector lattice lies in a free sublattice
on finitely many generators. We may therefore suppose that
\[
T:X\simeq \mathscr L(C)\to \operatorname{FVL}[\mathbb{R}^{k}]
\]
is injective. By duality, $T$ is induced by an $\ell$-map
$\varphi:\mathbb R^k\to C$. Its image is a closed polyhedral cone:
on each of finitely many closed convex polyhedral cones covering
$\mathbb R^k$, the map is linear. If $\varphi(\mathbb R^k)\neq C$,
we could choose $f\in\operatorname{FVL}[\mathbb R^n]$ with
$S(f)=\mathbb R^n\setminus\varphi(\mathbb R^k)$. Then $f|_C\neq0$, whereas
\[
T(f|_C)=f\varphi=0,
\]
a contradiction. Hence $C=\varphi(\mathbb R^k)$ is polyhedral.
\end{proof}

\begin{cor}\label{cor:locally-projective-freeprod}
\leanref{FreeProducts/AVLvsVL/Coproducts/Equivalences.lean}{225}{archimedean_freeProducts_iff_locally_projective}
Let $X$ be an Archimedean vector lattice. Then $X*_{\scriptscriptstyle{\mathcal{VL}}} Y$ is Archimedean for every Archimedean $Y$ if and only if every finitely generated sublattice of $X$ satisfies one of the equivalent conditions from Proposition \ref{prop:X-fg-projective}.
\end{cor}
\begin{proof}
Suppose first that $X*_{\scriptscriptstyle{\mathcal{VL}}}Y$ is
Archimedean for every Archimedean $Y$. If $A\subseteq X$ is finitely
generated and $B$ is finitely generated and Archimedean, then
Corollary~\ref{cor:freeprod-vl-injective} embeds
$A*_{\scriptscriptstyle{\mathcal{VL}}}B$ into
$X*_{\scriptscriptstyle{\mathcal{VL}}}B$. Hence the former is
Archimedean, and Proposition~\ref{prop:X-fg-projective} applies to $A$.

Conversely, suppose every finitely generated sublattice of $X$ is
projective, and let $Y$ be Archimedean. The free product
$X*_{\scriptscriptstyle{\mathcal{VL}}}Y$ is the directed union of the
embedded products $A*_{\scriptscriptstyle{\mathcal{VL}}}B$, where
$A\subseteq X$ and $B\subseteq Y$ are finitely generated. Each such
product is Archimedean by Proposition~\ref{prop:X-fg-projective}.
If $0\leq nz\leq w$ for all $n$, a single such product contains
$z,w$ and reflects these inequalities. Thus $z=0$, proving that the
full free product is Archimedean.
\end{proof}

\begin{ex}
Not every Archimedean vector lattice such that every finitely generated sublattice is projective embeds in a free vector lattice.
Let

\[
X=\{f\colon\mathbb N\to\mathbb R:f(\mathbb N)\text{ is finite}\},
\]

with pointwise vector-lattice operations. Then:

\begin{enumerate}
\def\labelenumi{\arabic{enumi}.}
\setlength{\itemsep}{0pt}\setlength{\parskip}{0pt}

\item \leanref{FreeProducts/Counterexamples/LocalProjectivity.lean}{121}{finiteRangeSequences_locally_projective}
  Every finitely generated vector sublattice of \(X\) embeds into a
  finitely generated free vector lattice;
\item \leanref{FreeProducts/Counterexamples/LocalProjectivity.lean}{243}{finiteRangeSequences_not_free_embeddable}
  \(X\) does not embed into \(\operatorname{FVL}[E]\) for any vector space $E$.
\end{enumerate}

Indeed, let \(f_1,\ldots,f_d\in X\). The map
\[
n\longmapsto(f_1(n),\ldots,f_d(n))
\]
has finite range. Its fibers give a finite partition
\(\mathbb N=P_1\sqcup\cdots\sqcup P_r\), and the vector sublattice
generated by \(f_1,\ldots,f_d\) is contained in the lattice of functions
constant on every \(P_i\). The latter lattice is isomorphic to
\(\mathbb R^r\), which is projective and therefore embeds into a
finitely generated free vector lattice. The same embedding restricts
to the sublattice generated by $f_1,\ldots,f_d$, proving the first assertion.

Now recall that a component of an
element \(e\geq0\) is an element \(w\) satisfying
\[
0\leq w\leq e
\quad\text{and}\quad
w\wedge(e-w)=0.
\]

For any set $S$, every positive element of $\operatorname{FVL}[\mathbb R^S]$ has only finitely many components.
Indeed, let $e\in\operatorname{FVL}[\mathbb R^S]_+$. It depends on a finite set
$T\subseteq S$, so it may be regarded as an element of $\operatorname{FVL}[\mathbb R^T]$.
In the standard representation of $FVL$ as piecewise-linear functions, the open polyhedral cone
\[
P=\{x\in\mathbb R^T:e(x)>0\}
\]

has only finitely many connected components. If \(w\) is a component of
\(e\), choose a finite set \(U\supseteq T\) such that
$w\in\operatorname{FVL}[\mathbb R^U]$.
On \(P\times\mathbb R^{U\setminus T}\), the continuous function \(w/e\)
takes values only in \(\{0,1\}\). It is therefore constant on each
connected component. Outside this set, \(w=0\) because
\(0\leq w\leq e\). Thus a component is determined by a subset of the
finitely many connected components of \(P\), and only finitely many
components are possible.

On the other hand, the constant function \(\mathbf 1\in X\) has the
distinct components
\[
\mathbf 1_M(n)=
\begin{cases}
1,&n\in M,\\
0,&n\notin M,
\end{cases}
\qquad(M\subseteq\mathbb N).
\]
An embedding preserves components and their distinctness. Hence no
embedding from \(X\) into any \(\operatorname{FVL}[\mathbb R^S]\) can exist.
\end{ex}

\begin{ex}\label{ex:non-archimedean-freeprod}
\leanref{FreeProducts/Counterexamples/FreeBanachComparison.lean}{567}{freeBanach_pair_not_archimedean}
The $\mathcal{VL}$ free product of free Banach lattices need not be Archimedean. \leanref{FreeProducts/Counterexamples/FreeBanachComparison.lean}{588}{freeBanach_comparison_not_injective}
Thus, for Banach spaces $E$ and $F$, the map $\mathcal{I}$ in the diagram
\[
\xymatrix{\operatorname{FVL}[E\oplus F]\simeq\operatorname{FVL}[E]*_{\scriptscriptstyle{\mathcal{VL}}} \operatorname{FVL}[F] \ar@{^(->}[d] \ar@{^(->}[r] & \operatorname{FBL}[E]*_{\scriptscriptstyle{\mathcal{BL}}} \operatorname{FBL}[F]\simeq\operatorname{FBL}[E\oplus_1 F]  \\ \operatorname{FBL}[E]*_{\scriptscriptstyle{\mathcal{VL}}} \operatorname{FBL}[F] \ar[ru]_{\mathcal{I}} & }
\]
is not injective in general.

Indeed, take the closed cone
\[
C=\{(x,y,z)\in\mathbb R^3:z=\sqrt{x^2+y^2}\}.
\]
It is not polyhedral: its intersection with the plane $x_3=1$ is a circle,
and a convex subset of that circle contains at most one point. As in the proof of Lemma \ref{lema:cone-bad-partner}, consider the cones
\[
\Omega_C\quad\text{and}\quad D_C=\mathbb R^4\backslash\Omega_C.
\]

It follows that
\[
J_{C\times D_C}\subsetneq \Van(C\times D_C).
\]

To embed the two factors into free Banach lattices, choose continuous surjections
\[
q_C:S^1\to C\cap S^2,\quad\text{and}\quad q_{D_C}:S^3\to D_C\cap S^3.
\]
For the first map, take $q_C(s,t)=2^{-1/2}(s,t,1)$.
Observe that $D_C\cap S^3$ is homeomorphic to the closed Euclidean
ball $B^3$: writing $a(u)=\operatorname{dist}_\infty(u,C)$ for $u\in S^2$, a homeomorphism is
\[
\Phi(\rho u)=
\frac{(u,a(u)+1-\rho^{-1})}
{\sqrt{1+(a(u)+1-\rho^{-1})^2}}
\quad(0<\rho\leq1),\qquad \Phi(0)=(0,-1).
\]
Every point of $D_C\cap S^3$ other than $(0,-1)$ has the unique form
$(u,s)/\sqrt{1+s^2}$ with $s\leq a(u)$, corresponding to
$\rho=(a(u)+1-s)^{-1}$. Boundedness of $a$ gives continuity at $\rho=0$,
so the continuous bijection is a homeomorphism.
The projection $P:S^3\to B^3$ onto the first three coordinates is surjective;
hence $q_{D_C}=\Phi P:S^3\to D_C\cap S^3$ is a continuous surjection.

Pullback along $q_C$ and $q_{D_C}$ gives embeddings
\[
\mathscr L(C)\hookrightarrow C(S^1)\simeq\operatorname{FBL}[\ell_2^2],
\qquad
\mathscr L(D_C)\hookrightarrow C(S^3)\simeq\operatorname{FBL}[\ell_2^4].
\]
By Corollary \ref{cor:freeprod-vl-injective}, these induce an embedding
\[
\mathscr L(C)*_{\scriptscriptstyle{\mathcal{VL}}}\mathscr L(D_C)
\hookrightarrow
\operatorname{FBL}[\ell_2^2]*_{\scriptscriptstyle{\mathcal{VL}}}
\operatorname{FBL}[\ell_2^4].
\]
Therefore the latter cannot be Archimedean.
\end{ex}

\begin{prop}\label{prop:ck-avl-freeprod}
\leanref{FreeProducts/Archimedean/Embeddings.lean}{141}{exists_join_generated_equiv}
Let $X_1$, $X_2$ be nonzero
Archimedean vector lattices with specified strong units $u_1,u_2$.
Identify their order-unit completions with $C(K_1)$ and $C(K_2)$,
respectively, so that $u_i$ corresponds to $1_{K_i}$. Then
\[
X_1*_{\scriptscriptstyle{\mathcal{AVL}}} X_2
\]
is lattice isomorphic to the sublattice
\[
\lat(\varphi_1(X_1),\varphi_2(X_2))\subset C(K_1*K_2).
\]
\end{prop}

\begin{proof}
The proof is analogous to \cite[Theorem 6.2]{freeprods_BL}.
\end{proof}

Although $\mathcal{AVL}$ lacks the amalgamation property, its free products
still preserve embeddings. The join representation gives a proof for finitely
generated sublattices, and the general case follows by assembling these
finite pieces.

\begin{prop} \label{prop:injectivity-avl-freeprod}
\leanref{FreeProducts/Archimedean/Embeddings.lean}{596}{induced_injective}
Let $Y_i\subset X_i$ be Archimedean vector lattices, for $i=1,2$. Then the induced lattice homomorphism
\[
Y_1*_{\scriptscriptstyle{\mathcal{AVL}}} Y_2\to X_1*_{\scriptscriptstyle{\mathcal{AVL}}} X_2
\]
is injective.
\end{prop}

\begin{proof}
First suppose that all four lattices are finitely generated.
Choose strong units in $Y_i$ and $X_i$. The inclusions are bounded for
the associated order-unit norms and extend to lattice homomorphisms
$J_i:C(K_i)\to C(L_i)$ between their completions. These extensions are
injective. Indeed, $Y_i$ is order dense in $C(K_i)$, since it is
uniformly dense and contains the unit. A nonzero kernel of $J_i$ would
therefore contain a nonzero positive element of $Y_i$, contradicting
injectivity of the original inclusion. Theorem~\ref{theo:ck-injectivity}
makes the induced map on the Banach lattice free products injective.
Restricting it to the sublattices in Proposition~\ref{prop:ck-avl-freeprod}
proves the assertion in this case. 

For arbitrary $X_1,X_2$, consider the directed system
$A_1*_{\scriptscriptstyle{\mathcal{AVL}}}A_2$, where $A_i$ ranges over
the finitely generated sublattices of $X_i$. Whenever $A_i\subseteq B_i$
for $i=1,2$, the transition map
$A_1*_{\scriptscriptstyle{\mathcal{AVL}}}A_2\to
B_1*_{\scriptscriptstyle{\mathcal{AVL}}}B_2$ is the lattice homomorphism
induced by the inclusions of the factors. These maps are compatible by
uniqueness in the universal property. The finite case shows that
all transition maps are embeddings. Identify these lattices along the
transition maps and take their directed union $V$. This union is
Archimedean: if $0\leq nz\leq w$ for every $n$, a single member of the directed system contains
$z,w$ and these inequalities hold there so $z=0$.

The canonical maps $X_i\to V$ are compatible with the maps at each
stage. For an Archimedean target $Z$ and lattice homomorphisms
$T_i:X_i\to Z$, the induced maps on the directed system agree on
overlaps, by uniqueness in their universal properties. They therefore
give a unique lattice homomorphism $V\to Z$ extending both $T_i$.
Thus $V$ is $X_1*_{\scriptscriptstyle{\mathcal{AVL}}}X_2$.

Apply the same description to $Y_1,Y_2$. Each finite generated sublattice of their
free product embeds into a finite generated sublattice for $X_1,X_2$, and hence into
$V$. Since every element belongs to a finite generated sublattice, the induced map
on the full free products is injective.
\end{proof}

\section{Consequences for the free product of Banach lattices}\label{sec:banach-free-products}

\subsection{Maximal norm and completion}
We now identify the algebraic sublattice generated by the factors inside
the Banach lattice free product. It is their Archimedean free product,
and the Banach lattice free product is recovered by completing it under
the maximal compatible lattice norm.

\begin{theo}\label{theo:avl-into-bl}
\leanref{FreeProducts/Consequences/NormAndCompletion.lean}{529}{archimedean_completion}
Given a family $(X_i)_{i\in I}$ of Banach lattices, the lattice homomorphism
\[
\mathcal I:=\overline{\bigfree} \operatorname{id}_{X_i}:\bigfree\nolimits_{\scriptscriptstyle{\mathcal{AVL}}} X_i\to \bigfree\nolimits_{\scriptscriptstyle{\mathcal{BL}}} X_i
\]
is injective. Let $\mathcal M$ be the set of lattice seminorms $p$ on $\bigfree\nolimits_{\scriptscriptstyle{\mathcal{AVL}}} X_i$ such that \[
p(x_i)\leq\|x_i\|_{X_i} \]
for all $x_i\in X_i$ and $i\in I$. If we endow
$\bigfree\nolimits_{\scriptscriptstyle{\mathcal{AVL}}} X_i$
with the lattice norm
\[
\|z\|=\sup_{p\in\mathcal M} p(z),
\]
then $\mathcal I$ is a lattice isometric embedding with dense range.
\end{theo}
\begin{proof}
We first prove injectivity when $I=\{1,2\}$, writing the factors as $X,Y$.

Let $0\neq z\in X*_{\scriptscriptstyle{\mathcal{AVL}}} Y$. Write
\[
z=\Phi\big(\varphi_X(x_1),\ldots,\varphi_X(x_n),\varphi_Y(y_1),\ldots,\varphi_Y(y_m)\big )
\]
and put
\[
u=|x_1|\vee\ldots\vee|x_n|,\qquad v=|y_1|\vee\ldots\vee|y_m|.
\]
Let $X_u$ and $Y_v$ be the principal ideals generated by $u$ and $v$, respectively, which become Banach lattices when equipped with the corresponding order unit norms. We have a commutative diagram
\[
\xymatrix{ X_u*_{\scriptscriptstyle{\mathcal{AVL}}} Y_v \ar[r] \ar[d] & X_u*_{\scriptscriptstyle{\mathcal{BL}}}Y_v \ar[d] \\ X*_{\scriptscriptstyle{\mathcal{AVL}}} Y \ar[r]^{\mathcal I} & X*_{\scriptscriptstyle{\mathcal{BL}}} Y}
\]
The element
\[
0\neq \tilde z=\Phi\big(\varphi_{X_u}(x_1),\ldots,\varphi_{X_u}(x_n),\varphi_{Y_v}(y_1),\ldots,\varphi_{Y_v}(y_m)\big)
\]
in $X_u*_{\scriptscriptstyle{\mathcal{AVL}}}Y_v$ maps to $z$ under
the left vertical arrow. The top horizontal arrow is injective by
Proposition~\ref{prop:ck-avl-freeprod}, and the right vertical arrow
is injective by Theorem~\ref{theo:ck-injectivity}, using the order-unit
representations of $X_u$ and $Y_v$.
Commutativity gives $\mathcal I(z)\neq0$.

The finite case follows by induction. Embed the Archimedean free
product of the first factors into their Banach lattice free product
using the induction hypothesis, apply
Proposition~\ref{prop:injectivity-avl-freeprod} to adjoin the last
factor, and then use the two-factor case. For arbitrary $I$, every
element involves finitely many factors. The canonical maps from a
finite subproduct into the full product have left inverses obtained
by setting the other factors to zero, in both categories. Hence the
finite case gives injectivity for arbitrary families.

The norm inherited from $\bigfree\nolimits_{\scriptscriptstyle{\mathcal{BL}}} X_i$ belongs to $\mathcal M$. Every $z\in \bigfree\nolimits_{\scriptscriptstyle{\mathcal{AVL}}} X_i$ can be dominated by a finite sum of elements from the free factors, so the supremum is finite and defines a lattice norm on $\bigfree\nolimits_{\scriptscriptstyle{\mathcal{AVL}}} X_i$. Let $Z_0$ be its norm completion. For every family of contractive Banach lattice homomorphisms $T_i:X_i\to Z$, the induced vector lattice homomorphism
\[
\bigfree T_i:\bigfree\nolimits_{\scriptscriptstyle{\mathcal{AVL}}} X_i \to Z
\]
defines an element of $\mathcal M$ via $p(z)=\|(\bigfree T_i)(z)\|_Z$,
so $\bigfree T_i$ is contractive and extends to a unique contractive lattice homomorphism from $Z_0$. This is the universal property of free products of Banach lattices, and \cite[Remark 3.2]{freeprods_BL} completes the proof.
\end{proof}

\subsection{Order density and regularity}
A sublattice $A$ of $X$ is order dense if every $0<x\in X$ dominates
some $0<a\in A$. It is regular if every net decreasing to zero in
$A$ also decreases to zero in $X$. Every order dense sublattice is
regular.

The completion approach in the previous subsection also admits a local description: every element comes from
a free product of principal ideals in the two factors. Since these ideals
have strong units, this brings the topological join approach into play for elements of
the completed free product as well.

\begin{prop}\label{prop:any-element-range-principal-ideals}
\leanref{FreeProducts/Consequences/PrincipalLocalization.lean}{453}{exists_mem_range_principal_freeProduct}
Let $X,Y$ be Banach lattices. For every $z\in X*_{\scriptscriptstyle{\mathcal{BL}}}Y$, there are elements $u\in X_+$, $v\in Y_+$ such that $z$ belongs to the range of the canonical map
\[
X_u*_{\scriptscriptstyle{\mathcal{BL}}}Y_v\to X*_{\scriptscriptstyle{\mathcal{BL}}}Y.
\]
Here the principal ideals $X_u$ and $Y_v$ carry their order-unit norms.
\end{prop}

\begin{proof}
We first establish a formula for the free-product norm on
$X*_{\scriptscriptstyle{\mathcal{AVL}}}Y$. For
$w\in X*_{\scriptscriptstyle{\mathcal{AVL}}}Y$, set
\[
p(w)=\inf\bigl\{\|a\|_X+\|b\|_Y:
a\in X_+,\ b\in Y_+,\ |w|\leq\varphi_X(a)+\varphi_Y(b)\bigr\}.
\]
Every element of the algebraic sublattice generated by the factors is
dominated in modulus by a sum of positive elements from the factors, so
$p(w)$ is finite. Adding admissible majorants gives subadditivity, and
rescaling them gives absolute homogeneity. Moreover, $p(w)\leq p(w')$
whenever $|w|\leq|w'|$. Thus $p$ is a lattice seminorm, and
\[
p(\varphi_X(a))\leq\|a\|_X,\qquad
p(\varphi_Y(b))\leq\|b\|_Y
\qquad(a\in X,\ b\in Y).
\]
Theorem~\ref{theo:avl-into-bl} therefore gives $p(w)\leq\|w\|_*$.
Conversely, every pair $a,b$ admissible in the definition of $p(w)$ satisfies
\[
\|w\|_*\leq\|\varphi_X(a)+\varphi_Y(b)\|_*
\leq\|a\|_X+\|b\|_Y.
\]
Taking the infimum yields $\|w\|_*\leq p(w)$, so $p(w)=\|w\|_*$.

Now fix $z\in X*_{\scriptscriptstyle{\mathcal{BL}}}Y$.
By Theorem~\ref{theo:avl-into-bl}, choose
\(z_n\in X*_{\scriptscriptstyle{\mathcal{AVL}}}Y\) with \[
\|z-z_n\|_*<4^{-n}
\qquad(n\geq1).
\] Set \(z_0=0\) and \[
d_n=z_n-z_{n-1}.
\] Then \[
\sum_{n=1}^\infty 2^n\|d_n\|_*<\infty,
\] because for \(n\geq2\), \[
\|d_n\|_*\leq\|z-z_n\|_*+\|z-z_{n-1}\|_*<4^{-n}+4^{-(n-1)}.
\] Also \[
z=\sum_{n=1}^{\infty}d_n
\] in \(X*_{\scriptscriptstyle{\mathcal{BL}}}Y\).

Since $p(d_n)=\|d_n\|_*$, for each \(n\) we can choose elements
$u_n\in X_+$, $v_n\in Y_+$ such that
$|d_n|\leq \varphi_X(u_n)+\varphi_Y(v_n)$ and \[
\|u_n\|_X+\|v_n\|_Y\leq \|d_n\|_*+4^{-n}.
\] Then \[
\sum_{n=1}^{\infty}2^n(\|u_n\|_X+\|v_n\|_Y)<\infty.
\]
List the elements of $X$ and $Y$ occurring in chosen lattice-linear
expressions for the $z_n$ as sequences $(x_n)$ and $(y_n)$, respectively.
Define
\[
x=\sum_{n=1}^\infty \frac{|x_n|}{2^n(1+\|x_n\|_X)},\qquad
y=\sum_{n=1}^\infty \frac{|y_n|}{2^n(1+\|y_n\|_Y)}.
\]
Both series converge in norm, since the norm of each summand is at most
$2^{-n}$. Moreover,
\[
|x_n|\leq 2^n(1+\|x_n\|_X)x,\qquad
|y_n|\leq 2^n(1+\|y_n\|_Y)y.
\]

Let \[
u:=x\vee\sum_{n=1}^{\infty}2^n u_n\in X_+,\qquad v:=y\vee\sum_{n=1}^\infty 2^n v_n.
\]
In particular, $u_n\leq 2^{-n}u$, $v_n\leq 2^{-n}v$, and every
$x_n$ belongs to $X_u$, while every $y_n$ belongs to $Y_v$.

The principal ideals $X_u$ and $Y_v$, with their order-unit norms
$\|a\|_u=\inf\{\lambda\geq0:|a|\leq\lambda u\}$ and
$\|b\|_v=\inf\{\lambda\geq0:|b|\leq\lambda v\}$, are Banach lattices.
Their inclusions into $X$ and $Y$ are bounded, since
$\|a\|_X\leq\|u\|_X\|a\|_u$ and
$\|b\|_Y\leq\|v\|_Y\|b\|_v$. Consider the induced bounded lattice
homomorphism
\[
j:X_u*_{\scriptscriptstyle{\mathcal{BL}}}Y_v\to X*_{\scriptscriptstyle{\mathcal{BL}}}Y.
\]
The order-unit representations of $X_u$ and $Y_v$ as
$C(K)$ spaces allow us to apply Theorem~\ref{theo:ck-injectivity}, so
$j$ is injective.

Every \(d_n\) has a preimage \[
\widetilde d_n\in X_u*_{\scriptscriptstyle{\mathcal{AVL}}} Y_v\subset X_u*_{\scriptscriptstyle{\mathcal{BL}}}Y_v,
\]
obtained by taking the difference of the chosen expressions for $z_n$ and
$z_{n-1}$ in these principal ideals. Since an injective lattice homomorphism
reflects order, the bound on \(d_n\) gives \[
|\widetilde d_n|
\leq
\varphi_{X_u}(u_n)+\varphi_{Y_v}(v_n),
\] so the free-product norm satisfies \[
\|\widetilde d_n\|_{X_u*_{\scriptscriptstyle{\mathcal{BL}}}Y_v}
\leq
\|u_n\|_u+\|v_n\|_v
\leq
2^{-n+1}.
\] Thus \(\sum_n\widetilde d_n\) converges in
$X_u*_{\scriptscriptstyle{\mathcal{BL}}}Y_v$. Define \[
\widetilde z=\sum_{n=1}^{\infty}\widetilde d_n.
\]

By continuity of $j$, \[
j\widetilde z=\sum_{n=1}^{\infty}d_n=z.
\]
\end{proof}

\begin{cor}\label{cor:order-denseness}
\leanref{FreeProducts/Consequences/Families.lean}{234}{archimedean_induced_orderDense}
Let $X',Y'$ be order dense sublattices of the Banach lattices $X$ and $Y$. Then
\[
X'*_{\scriptscriptstyle{\mathcal{AVL}}} Y'\subset X*_{\scriptscriptstyle{\mathcal{BL}}}Y
\]
is an order dense sublattice. \leanref{FreeProducts/Consequences/Families.lean}{250}{archimedean_induced_regular}
In particular, it is regular.
The same conclusion holds for any finite family of factors and order
dense sublattices.
\end{cor}

\begin{proof}
Let $0<z\in X*_{\scriptscriptstyle{\mathcal{BL}}}Y$. By Proposition \ref{prop:any-element-range-principal-ideals}, choose $u$ and $v$ such that $z$ belongs to the image of $j:X_u*_{\scriptscriptstyle{\mathcal{BL}}}Y_v\to X*_{\scriptscriptstyle{\mathcal{BL}}}Y$. Let $\tilde z$ be the preimage of $z$; then $\tilde z>0$. By \cite[Proposition 6.7]{freeprods_BL}, we can choose \[
w\in (X'\cap X_u)*_{\scriptscriptstyle{\mathcal{AVL}}}(Y'\cap Y_v)\subset X_u*_{\scriptscriptstyle{\mathcal{AVL}}} Y_v
\]
such that $0<w\leq \tilde z$. Injectivity of $j$ then implies
\[
0<j w\leq z, \quad j w\in X'*_{\scriptscriptstyle{\mathcal{AVL}}} Y'.
\]
For finitely many factors, argue by induction: the Archimedean free
product of the first factors is order dense in their Banach lattice
free product, so the two-factor assertion applies when the last
factor is adjoined. Associativity identifies the resulting sublattice
with the full Archimedean free product.
\end{proof}

\begin{ex}
\leanref{FreeProducts/Counterexamples/OrderDensity.lean}{190}{countable_freeProduct_not_orderDense}
Order density need not hold for infinitely many free factors. \leanref{FreeProducts/Counterexamples/OrderDensity.lean}{167}{countable_freeProduct_closed_ideal_disjoint}
In fact, even the weaker property that closed principal ideals contain nonzero elements from the $\mathcal{AVL}$ free product can fail. Take $X_n=\operatorname{FBL}[\mathbb R]$ for every $n\in \mathbb N$. Then
\[
\bigfree\nolimits_{\scriptscriptstyle{\mathcal{BL}}} X_n\simeq \operatorname{FBL}[\ell_1]
\]
as in \cite[Proposition 5.5]{freeprods_BL}. Let
\[
x=\sum_{n=1}^\infty \frac{e_n}{2^n}\in\ell_1\quad\text{and}\quad h=|\delta_x|\in \operatorname{FBL}[\ell_1].
\]
The closed ideal generated by $h$, denoted by $\overline{I_h}$, is contained in the ideal of functions vanishing on
\[
Z(x)=\{y^*\in\ell_\infty:y^*(x)=0\}.
\]
Now let \[
f\in \overline{I_h}\cap \bigfree\nolimits_{\scriptscriptstyle{\mathcal{VL}}} X_n.
\]
There is a finite set $S\subset\mathbb N$ such that $f$ depends only on the coordinates $(y^*(e_n))_{n\in S}$. We show that $f=0$. Fix
arbitrary scalar values \((\alpha_n)_{n\in S}\). Choose \(m\notin S\)
and define \(y^*\in\ell_\infty\) by \[
y^*(e_n)=\alpha_n\quad(n\in S),
\qquad
y^*(e_m)=-2^m\sum_{n\in S}2^{-n}\alpha_n,
\] and \(y^*(e_k)=0\) for all other \(k\). Then \(y^*\) is bounded and
\[
y^*(x)=\sum_{n\in S}2^{-n}\alpha_n
+2^{-m}y^*(e_m)=0,
\] so \(y^*\in Z(x)\). Since \(f\in \overline{I_h}\), it follows that \(f(y^*)=0\).
But \(f\) depends only on the coordinates in \(S\), and the values
\((\alpha_n)_{n\in S}\) were arbitrary. Therefore \(f\) vanishes for
every possible value of its finite coordinate variables, so \(f=0\).
\end{ex}

\begin{cor}\label{cor:free-factor-regular}
\leanref{FreeProducts/Consequences/Families.lean}{274}{inclusion_regular}
Let $(X_i)_{i\in I}$ be a family of Banach lattices and $i_0\in I$. Then
\[
X_{i_0}\subset \bigfree\nolimits_{\scriptscriptstyle{\mathcal{BL}}} X_i
\]
is a regular sublattice.
\end{cor}

\begin{proof} 
Since
\[
\bigfree\nolimits_{\scriptscriptstyle{\mathcal{BL}}} X_i\simeq X_{i_0}*_{\scriptscriptstyle{\mathcal{BL}}}\underset{\{i\neq i_0\}}{\bigfree\nolimits_{\scriptscriptstyle{\mathcal{BL}}}} X_i,
\]
the claim follows from Corollary \ref{cor:order-denseness} and \cite[Proposition 6.10]{freeprods_BL}.
\end{proof}

\subsection{Preservation of injectivity}

\begin{prop}\label{prop:injectivity-bl-freeprod}
\leanref{FreeProducts/Consequences/Families.lean}{430}{induced_injective}
Let $(T_i:X_i\to Y_i)_{i\in I}$ be a uniformly bounded family of injective lattice homomorphisms between Banach lattices. Then
\[
\bigfree T_i:\bigfree\nolimits_{\scriptscriptstyle{\mathcal{BL}}} X_i\to \bigfree\nolimits_{\scriptscriptstyle{\mathcal{BL}}} Y_i
\]
is injective.
\end{prop}
\begin{proof}
For $F\subseteq I$, write
\[
X_F=\underset{i\in F}{\bigfree\nolimits_{\scriptscriptstyle{\mathcal{BL}}}}X_i,
\qquad
Y_F=\underset{i\in F}{\bigfree\nolimits_{\scriptscriptstyle{\mathcal{BL}}}}Y_i,
\]
and let $T_F:X_F\to Y_F$ be the homomorphism induced by $(T_i)_{i\in F}$.
The uniform bound on the $T_i$ ensures that these maps are bounded.

First, suppose that $I=\{1,2\}$. By
Theorem~\ref{theo:avl-into-bl} and
Proposition~\ref{prop:injectivity-avl-freeprod}, the restriction of $T_I$ to
$X_1*_{\scriptscriptstyle{\mathcal{AVL}}}X_2$ is injective. This sublattice
is order dense in $X_I$ by Corollary~\ref{cor:order-denseness}. If
$0\neq z\in\ker T_I$, choose
$w\in X_1*_{\scriptscriptstyle{\mathcal{AVL}}}X_2$ with $0<w\leq|z|$.
Then
\[
0\leq T_I(w)\leq T_I(|z|)=|T_I(z)|=0,
\]
contradicting injectivity on that sublattice. The general finite case follows by induction using associativity of Banach lattice free products.

Now let $I$ be infinite. For every finite $F\subseteq I$, the universal
property gives a canonical contractive map $j_F^X:X_F\to X_I$ and a
contractive lattice homomorphism $P_F^X:X_I\to X_F$ induced by the
canonical maps on the factors indexed by $F$ and the zero maps on the
remaining factors. Since $P_F^Xj_F^X=\operatorname{id}_{X_F}$, the map $j_F^X$ is an
isometric lattice embedding, and we identify $X_F$ with its image in
$X_I$. Define $P_F^Y:Y_I\to Y_F$ in the same way.

We claim that if $z\in X_I$ satisfies $P_F^Xz=0$ for every finite
$F\subseteq I$, then $z=0$. By Theorem~\ref{theo:avl-into-bl}, choose
$z_n\in\underset{i\in I}{\bigfree\nolimits_{\scriptscriptstyle{\mathcal{AVL}}}}X_i$ with
\[
\|z-z_n\|\leq 1/n
\]
and let $F_n$ be a finite set of indices containing all factors occurring
in a lattice-linear expression for $z_n$. Then $z_n\in X_{F_n}$ and
$P_{F_n}^Xz_n=z_n$. Our assumption gives $P_{F_n}^Xz=0$, so
\[
\|z_n\|=\|P_{F_n}^X(z_n-z)\|\leq\|z_n-z\|\leq 1/n.
\]
Consequently, $\|z\|\leq\|z-z_n\|+\|z_n\|\leq2/n$ for every $n$,
which proves the claim.

For every finite $F\subseteq I$, the diagram
\[
\xymatrix{
X_I \ar[r]^{T_I} \ar[d]_{P_F^X} & Y_I \ar[d]^{P_F^Y} \\
X_F \ar[r]_{T_F} & Y_F
}
\]
commutes, since both compositions agree on every free factor. If
$z\in\ker T_I$, it follows that
\[
T_F(P_F^Xz)=P_F^Y(T_Iz)=0.
\]
The finite case shows that $T_F$ is injective, so $P_F^Xz=0$ for every
finite $F$. The claim now gives $z=0$, proving that $T_I$ is injective.
\end{proof}

\section{Vector lattices with countable spectra}
\label{sec:countable-spectrum}

In this section, we apply the amalgamation results to extension problems for vector lattices with countably many prime ideals. These problems are analogous to those concerning
injectivity and universal disposition for separable Banach lattices studied in
\cite{amalg}, but here size is measured by the cardinality of the prime
spectrum. All vector lattices are real, and no Archimedean assumption is
imposed unless explicitly stated.
For background on prime ideals in vector lattices, see \cite[Chapter~5]{luxemburg} or \cite{JohnsonKist1962}.

Throughout this section, the spectrum of a vector lattice $X$ consists of its \emph{proper}
prime ideals:
\[
 \operatorname{Spec}_p(X)=
 \{P\subsetneq X:P\text{ is an ideal and }X/P\text{ is totally ordered}\}.
\]
Let $\mathcal C$ denote the class of vector lattices with countable spectrum:
\[
 \mathcal C=\{X:|\operatorname{Spec}_p(X)|\leq\aleph_0\}.
\]
Thus countable spectrum includes finite spectra, and
$\operatorname{Spec}_p(0)=\varnothing$. 

\subsection{The size of a vector lattice with countable spectrum}\label{subsec:countable spectrum}

For a subset $M$ of a vector lattice $X$, write $I_X(M)$ for the ideal
generated by $M$.

\begin{lema}\label{cpi:hereditary}
Let $Y\subset X$ be a vector sublattice.
\begin{enumerate}
\item \leanref{FreeProducts/CountableSpectra/Basic.lean}{500}{mem_generated_ideal_iff}
For every ideal $J$ of $Y$,
\[
 I_X(J)\cap Y=J.
\]
\item \leanref{FreeProducts/CountableSpectra/Basic.lean}{207}{spectrum_cardinal_sublattice}
\[
 |\operatorname{Spec}_p(Y)|\leq|\operatorname{Spec}_p(X)|.
\]
\item \leanref{FreeProducts/CountableSpectra/Basic.lean}{172}{spectrum_cardinal_le_of_surjective}
If $K$ is an ideal of $X$, then
\[
 |\operatorname{Spec}_p(X/K)|\leq|\operatorname{Spec}_p(X)|.
\]
\end{enumerate}
In particular, $\mathcal C$ is closed under sublattices, quotients and
homomorphic images.
\end{lema}
\begin{proof}
Since $J$ is an ideal of $Y$,
\[
 I_X(J)=\{x\in X:|x|\leq j\text{ for some }j\in J_+\}.
\]
If $x$ also belongs to $Y$, solidity in $Y$ gives $x\in J$.
This proves the first assertion and yields the canonical embedding
$Y/J\to X/I_X(J)$.

By Proposition~\ref{prop:prime-ideals-extcon}, every proper prime ideal of $Y$
extends to a prime ideal of $X$. Such an extension is proper, and distinct
prime ideals of $Y$ require distinct extensions. Notice that restriction in
the other direction can give the whole sublattice: for
$P\in\operatorname{Spec}_p(X)$, either $P\cap Y=Y$ or
$P\cap Y\in\operatorname{Spec}_p(Y)$.
Proper prime ideals of $X/K$ correspond bijectively to proper prime ideals
of $X$ containing $K$. Finally, homomorphic images are quotients by kernels.
\end{proof}

We also recall prime ideal separation. If $0\neq x\in X$, extend the zero
prime ideal of the one-dimensional sublattice $\mathbb R|x|$ to a prime ideal of
$X$ which does not contain $x$. Thus, the natural homomorphism
\begin{equation}\label{eq:spec-injective}
 \begin{aligned}
 \iota:&X&\longrightarrow&\prod_{P\in\operatorname{Spec}_p(X)}X/P,\\
 &x&\longmapsto&(x+P)_{P\in\operatorname{Spec}_p(X)}
 \end{aligned}
\end{equation}
is injective, so $X$ can be regarded as a sublattice of a product of totally ordered vector lattices.

There is no distinction between the category of totally ordered vector spaces and that of totally ordered vector lattices. For a totally ordered real vector space $L$, denote its set of
\emph{Archimedean classes} by
\[
\Gamma_L=(L\setminus\{0\})/{\asymp},
\]
where
\[
 x\asymp y
 \quad\Longleftrightarrow\quad
 |x|\leq C|y|\ \text{and}\ |y|\leq D|x|
 \quad\text{for some }C,D>0.
\]
Write $v(x)$ for the class of $x$, and order the classes by
\[
 v(x)\geq v(y)
 \quad\Longleftrightarrow\quad
 |x|\leq C|y|\text{ for some }C>0,
\]
and set $v(0)=\infty$. For nonzero $x,y$,
\[
 v(x)>v(y)
 \quad\Longleftrightarrow\quad
 |x|<\epsilon|y|\text{ for every }\epsilon>0.
\]
We have $v(rx)=v(x)$ for $r\neq0$ and
$v(x+y)\geq\min\{v(x),v(y)\}$. If the two values are distinct, say $v(x)<v(y)$, then $v(x+y)=v(x)$, and $x+y>0$ exactly when $x>0$.

\begin{lema}\label{cpi:chain-cardinality}
For every nonzero totally ordered real vector space $L$,
\begin{enumerate}
\item \leanref{FreeProducts/CountableSpectra/Basic.lean}{306}{archimedeanClasses_cardinal_le_spectrum}
$|\Gamma_L|\leq|\operatorname{Spec}_p(L)|$;
\item \leanref{FreeProducts/CountableSpectra/Basic.lean}{415}{cardinal_le_power_archimedeanClasses}
$|L|\leq\mathfrak c^{|\Gamma_L|}$.
\end{enumerate}
\end{lema}

\begin{proof}
Choose $e_\gamma>0$ in each class $\gamma\in\Gamma_L$, and let
\[
 D_\gamma=\{x:|x|\leq r e_\gamma\text{ for some }r>0\},
 \qquad
 N_\gamma=\{x:n|x|\leq e_\gamma\text{ for every }n\geq1\}.
\]
These are ideals, and $N_\gamma$ is proper and prime in $L$.
If $\gamma>\eta$, then
$e_\gamma\in N_\eta\setminus N_\gamma$, proving the first inequality.

For $x\in D_\gamma$ define
\[
 q_\gamma(x)=\sup\{r\in\mathbb R:r e_\gamma\leq x\}.
\]
This is a well-defined real number, and, for every $\varepsilon>0$,
\[
 (q_\gamma(x)-\varepsilon)e_\gamma\leq x
       \leq(q_\gamma(x)+\varepsilon)e_\gamma.
\]
By definition, if $a,b\in\mathbb R$ and $z\in D_\gamma$ satisfy
$ae_\gamma\leq z\leq be_\gamma$, then $a\leq q_\gamma(z)\leq b$.
For $x,y\in D_\gamma$, adding the preceding bounds gives
\[
 (q_\gamma(x)+q_\gamma(y)-2\varepsilon)e_\gamma\leq x+y
       \leq(q_\gamma(x)+q_\gamma(y)+2\varepsilon)e_\gamma.
\]
Consequently,
\[
 |q_\gamma(x+y)-q_\gamma(x)-q_\gamma(y)|\leq2\varepsilon.
\]
Since this holds for every $\varepsilon>0$,
we obtain $q_\gamma(x+y)=q_\gamma(x)+q_\gamma(y)$. Similarly, multiplying the bounds for $x$
by $r\in\mathbb R$, reversing the inequalities when $r<0$, yields
\[
 (rq_\gamma(x)-|r|\varepsilon)e_\gamma\leq rx
       \leq(rq_\gamma(x)+|r|\varepsilon)e_\gamma.
\]
Thus $|q_\gamma(rx)-rq_\gamma(x)|\leq|r|\varepsilon$ for every $\varepsilon>0$,
which implies $q_\gamma(rx)=rq_\gamma(x)$. This proves that $q_\gamma$ is linear.

If $q_\gamma(x)=0$, the bounds for $x$ give $|x|\leq\varepsilon e_\gamma$ for
every $\varepsilon>0$. Taking $\varepsilon=1/n$ shows that
$n|x|\leq e_\gamma$ for every $n\geq1$, so $x\in N_\gamma$. Conversely,
if $x\in N_\gamma$, then
\[
 -\frac1n e_\gamma\leq x\leq\frac1n e_\gamma\qquad(n\geq1),
\]
and hence $|q_\gamma(x)|\leq1/n$ for every $n\geq1$. Therefore $q_\gamma(x)=0$,
and $\ker q_\gamma=N_\gamma$.

Let $f_\gamma: L\to\mathbb R$ be a linear functional extending $q_\gamma$.
For $x\neq0$, taking $\gamma=v(x)$ gives
$x\in D_\gamma\setminus N_\gamma$, so $f_\gamma(x)\neq0$.
Consequently $x\mapsto(f_\gamma(x))_{\gamma\in\Gamma_L}$ is a linear
injection into $\mathbb R^{\Gamma_L}$, proving the second inequality.
\end{proof}

\begin{prop}\label{cpi:cardinality}
\leanref{FreeProducts/CountableSpectra/Basic.lean}{485}{cardinal_le_continuum}
Every vector lattice in $\mathcal C$ has cardinality at most
$\mathfrak c$. \leanref{FreeProducts/CountableSpectra/Existence.lean}{110}{exists_small_sublattice_representatives}
There is a set of at most $2^{\mathfrak c}$
representatives of the isomorphism types in $\mathcal C$, and of the
embedding types $A\hookrightarrow B$ with $B\in\mathcal C$.
\end{prop}
\begin{proof}
If $X\in\mathcal C$, each prime quotient $X/P$ has countable spectrum
by Lemma~\ref{cpi:hereditary}, and therefore has cardinality at most
$\mathfrak c^{\aleph_0}=\mathfrak c$ by
Lemma~\ref{cpi:chain-cardinality}. The embedding of $X$ into a
countable product of these quotients given in \eqref{eq:spec-injective} yields $|X|\leq\mathfrak c$.

Structures of cardinality at most $\mathfrak c$ can be coded on subsets
of a fixed set of size $\mathfrak c$. The vector-lattice language has
$\mathfrak c$ operation symbols, all of finite arity, so there are at
most $2^{\mathfrak c}$ such structures. Specifying a sublattice or an
embedding does not increase this bound.
\end{proof}

The cardinality bound in Proposition~\ref{cpi:cardinality} is sharp:
a vector lattice with countable spectrum can have Hamel dimension
$\mathfrak c$.

\begin{ex}\label{cpi:lex-product}
Countable spectrum does not imply countable generation. Consider
$L=\mathbb R^{\mathbb N}$ with the lexicographic order. Its proper ideals are $0$ and
\[
 L_n=\{x\in \mathbb R^{\mathbb N}:x_1=\cdots=x_n=0\}\qquad(n\geq1).
\]
Indeed, let $J$ be a nonzero proper ideal, and let $i$ be the least
index that occurs as the first nonzero coordinate of an element of
$J$. A vector in $J$ with first nonzero coordinate $i$ dominates,
up to a positive scalar, every vector whose first $i-1$ coordinates
vanish. Hence $J=L_{i-1}$, with $L_0=L$; properness forces $i\geq2$. \leanref{FreeProducts/Counterexamples/LexicographicSpectrum.lean}{84}{lex_sequences_countable_spectrum}
Thus
$L\in\mathcal C$. However, the vectors $(1,t,t^2,\ldots)$, $t>0$,
form a linearly independent family, so
$\dim_{\mathbb R}L=\mathfrak c$. Every linear subspace of a totally
ordered space is a sublattice. \leanref{FreeProducts/Counterexamples/LexicographicSpectrum.lean}{94}{lex_sequences_not_countably_generated}
Hence $L$ needs $\mathfrak c$ lattice
generators.
\end{ex}

\subsection{Universal disposition and the push-out construction}

\begin{defin}\label{cpi:injectivity}
\leanref{FreeProducts/CountableSpectra/Basic.lean}{140}{IsCountablePrimeInjective}
A vector lattice $Z$ is \emph{injective relative to $\mathcal C$}, or
\emph{countable-prime-injective}, if every lattice homomorphism
$T:A\to Z$, where $A$ is a sublattice of $B\in\mathcal C$, extends to a
lattice homomorphism $\widetilde T:B\to Z$.
\end{defin}

\begin{defin}\label{cpi:universal-disposition}
\leanref{FreeProducts/CountableSpectra/Basic.lean}{149}{HasCountablePrimeUniversalDisposition}
A vector lattice $U$ is \emph{of universal disposition for
$\mathcal C$} if, for every embedding $i:A\to B$ with $B\in\mathcal C$
and every embedding $j:A\to U$, there is an embedding $k:B\to U$
such that $ki=j$.
\end{defin}

The bound on the size of the prime spectrum is imposed on the ambient lattice $B$; its sublattice $A$ then also belongs to $\mathcal C$. Neither $Z$ nor $U$ is
required to have countable spectrum. Taking $A=0$ in Definition~\ref{cpi:universal-disposition}
shows that $U$ contains every member of $\mathcal C$ as a sublattice.

\begin{theo}\label{cpi:ud-injective}
\leanref{FreeProducts/CountableSpectra/InjectiveProperties.lean}{97}{isCountablePrimeInjective_of_universalDisposition}
Every vector lattice $U$ of universal disposition for $\mathcal C$ is
countable-prime-injective. \leanref{FreeProducts/CountableSpectra/InjectiveProperties.lean}{83}{extension_with_generated_kernel}
More precisely, for $A\subset B\in\mathcal C$
and a lattice homomorphism $T:A\to U$, there is an extension $\widetilde T:B\to U$ satisfying
\[
 \ker\widetilde T=I_B(\ker T).
\]
\end{theo}
\begin{proof}
Put $K=\ker T$ and $J=I_B(K)$. Lemma~\ref{cpi:hereditary} gives an
embedding $A/K\to B/J$, and $B/J\in\mathcal C$. The map induced by
$T$ on $A/K$ is an embedding into $U$, so universal disposition extends
it to an embedding $H:B/J\to U$. Composing $H$ with the quotient map
$B\to B/J$ gives the required extension, with kernel exactly $J$.
\end{proof}

\begin{ex}\label{cpi:closure-failures}
\leanref{FreeProducts/Counterexamples/FreeProductSpectrum.lean}{40}{real_pair_uncountable_spectrum}
The class $\mathcal C$ is not closed under free products, even with two
one-dimensional factors. Let $x,y$ be the canonical positive generators
of $\mathbb R*_{\scriptscriptstyle{\mathcal{VL}}}\mathbb R$. For each
$t>0$, the homomorphism $\phi_t$ into $\mathbb R$ determined by
$\phi_t(x)=1$ and $\phi_t(y)=t$ is onto and its kernel is a proper prime ideal.
These kernels are distinct, since $y-tx$ belongs to $\ker\phi_t$ but
not to $\ker\phi_s$ when $s\neq t$. Thus this free product has at least
$\mathfrak c$ proper prime ideals. 

\leanref{FreeProducts/Counterexamples/IncreasingUnionSpectrum.lean}{81}{cylinder_union_uncountable_spectrum}
Similarly, the class $\mathcal C$ is not closed under countable increasing unions. On
$K=\{0,1\}^{\mathbb N}$, let $X_n$ be the space of functions depending on the
first $n$ coordinates. Then $X_n\cong\mathbb R^{2^n}$ has finite
spectrum, whereas distinct points of $K$ give distinct evaluation
kernels on $\bigcup_n X_n$. This union has at least $\mathfrak c$
proper prime ideals.
\end{ex}

The quotient argument in Theorem~\ref{cpi:ud-injective} avoids any
assumption that a push-out remains in $\mathcal C$. For the existence result, we
use the bound on underlying cardinality instead.

\begin{theo}\label{cpi:ud-existence}
\leanref{FreeProducts/CountableSpectra/Existence.lean}{316}{exists_countablePrimeUniversalDisposition}
Every vector lattice $X$ embeds into a vector lattice $U$ of universal
disposition for $\mathcal C$. With
\[
 \lambda=\bigl(\max\{|X|,2^{\mathfrak c}\}\bigr)^{\mathfrak c},
\]
one can arrange that $|U|\leq\lambda$. \leanref{FreeProducts/CountableSpectra/Existence.lean}{291}{exists_small_universalDisposition}
In fact, $U$ can be chosen of
universal disposition for all vector lattices of underlying cardinality
at most $\mathfrak c$.
\end{theo}
\begin{proof}
We have $\lambda^{\mathfrak c}=\lambda$ and
$\mathfrak c^+\leq2^{\mathfrak c}\leq\lambda$. By the coding argument
in Proposition~\ref{cpi:cardinality}, choose at most $2^{\mathfrak c}$
representatives of embeddings $i:A\hookrightarrow B$ with
$|B|\leq\mathfrak c$.

Construct an increasing continuous chain
$(U_\alpha)_{\alpha<\mathfrak c^+}$, beginning with $U_0=X$, with all
stages of cardinality at most $\lambda$. At each successor stage, list all
pairs
\[
 (i:A\hookrightarrow B,\ j:A\hookrightarrow U_\alpha),
\]
where $i$ is one of the chosen representatives. There are at most
$2^{\mathfrak c}\lambda^{\mathfrak c}=\lambda$ such pairs. Process
them successively, forming the push-out of $i$ and the composite of
$j$ with the inclusion of $U_\alpha$ into the current lattice.
Proposition~\ref{prop:pushout-vl}, applied in both directions, makes
both canonical maps embeddings. At internal limit stages take unions.

A push-out is generated by the images of its two factors, so for
factors $V$ and $B$ its cardinality is at most
$\max\{\mathfrak c,|V|,|B|\}$. Therefore the whole successor stage
$U_{\alpha+1}$ has cardinality at most $\lambda$. At outer limits take
$U_\alpha=\bigcup_{\beta<\alpha}U_\beta$, preserving the same bound.

Let $U=\bigcup_{\alpha<\mathfrak c^+}U_\alpha$. For an embedding
$i:A\to B$ with $|B|\leq\mathfrak c$ and an embedding $j:A\to U$,
regularity of $\mathfrak c^+$ places $j(A)$ inside some $U_\alpha$.
Up to the chosen representative, this diagram was processed at the
next stage, giving the required embedding $B\to U$. Finally,
Proposition~\ref{cpi:cardinality} includes every member of $\mathcal C$
among the considered objects.
\end{proof}

\begin{cor}\label{cpi:ud-size}
\leanref{FreeProducts/CountableSpectra/Existence.lean}{348}{exists_nontrivial_countablePrimeUniversalDisposition}
There is a nonzero vector lattice $U$ of universal disposition for
$\mathcal C$ such that
\[
 |U|\leq2^{\mathfrak c},\qquad
 |\operatorname{Spec}_p(U)|\leq2^{2^{\mathfrak c}}.
\]
\leanref{FreeProducts/CountableSpectra/Existence.lean}{332}{exists_countablePrimeInjective}
Every vector lattice embeds into a countable-prime-injective vector
lattice.
\end{cor}
\begin{proof}
For the first assertion take $X=\mathbb R$ in
Theorem~\ref{cpi:ud-existence} and use
$|\operatorname{Spec}_p(U)|\leq2^{|U|}$. The second follows from that
result and Theorem~\ref{cpi:ud-injective}.
\end{proof}

The construction uses transfinite induction in $\mathfrak c^+$ to capture every
domain of cardinality at most $\mathfrak c$. Example~\ref{cpi:lex-product}
shows why countable generation cannot be used to replace this argument
by an inductive construction in $\omega_1$. 

\subsection{Ordered ultrapowers and relative injectivity}

For a nonprincipal ultrafilter $\mathcal U$ on $\mathbb N$, consider
\[
 X_{\mathcal U}=\mathbb R^{\mathbb N}/\mathcal U.
\]
Two sequences represent the same element when they agree on a set in
$\mathcal U$, and
\[
 [(x_n)]\leq[(y_n)]
 \quad\Longleftrightarrow\quad
 \{n:x_n\leq y_n\}\in\mathcal U.
\]
With coordinatewise operations and real scalars identified with
constant sequences, this is a totally ordered real vector space of
cardinality $\mathfrak c$. We will prove that it is
$\mathcal C$-injective. In particular, the target can be much
smaller than the bound supplied by the earlier transfinite push-out construction.

For subsets $L,R$ of a totally ordered space, write $L<R$ when
$\ell<r$ for every $\ell\in L$ and $r\in R$. We say that a vector lattice has the
\emph{countable interpolation property} if, whenever $L,R$ are countable
and $L<R$, there is $z$ with $L<z<R$. 

The following lemma is a consequence of the countable saturation of
ultrapowers over nonprincipal ultrafilters on $\mathbb N$
\cite[Theorem~5.6]{Keisler2010}, applied to the underlying order in the
language $\{<\}$. We include a direct proof for completeness.

\begin{lema}\label{cpi:interpolation}
\leanref{FreeProducts/CountableSpectra/Ultrapowers.lean}{161}{orderedUltrapower_hasCountableInterpolation}
$X_{\mathcal U}$ has the countable interpolation property.
\end{lema}

\begin{proof}
Let $L,R\subseteq X_{\mathcal U}$ be countable sets with $L<R$.
Suppose first that both sets are nonempty. Write
$L=\{\ell_j:j\ge1\}$ and $R=\{r_j:j\ge1\}$, repeating elements if
either set is finite, and choose representatives
\[
    \ell_j=[(\ell_{j,n})_n],\qquad r_j=[(r_{j,n})_n].
\]
Since $\ell_j<r_k$ for every $j,k$, the definition of the ultrapower
order gives
\[
    \{n:\ell_{j,n}<r_{k,n}\}\in\mathcal U.
\]
For $m,n\ge1$, put
\[
    a_m(n)=\max_{1\le j\le m}\ell_{j,n},\qquad
    b_m(n)=\min_{1\le j\le m}r_{j,n},
\]
and define
\[
    E_m=\{n\ge m:a_m(n)<b_m(n)\}
       =\{n:n\ge m\}\cap
         \bigcap_{1\le j,k\le m}\{n:\ell_{j,n}<r_{k,n}\}.
\]
Each $E_m$ belongs to $\mathcal U$, since $\mathcal U$ is closed under
finite intersections and contains every cofinite subset of $\mathbb N$.
Moreover, $E_{m+1}\subseteq E_m$.

For every $n\ge1$, let
\[
    k(n)=\max\bigl(\{0\}\cup
            \{m\in\{1,\ldots,n\}:n\in E_m\}\bigr),
\]
and set
\[
    z_n=\begin{cases}
      \dfrac{a_{k(n)}(n)+b_{k(n)}(n)}{2},& k(n)>0,\\[4pt]
      0,& k(n)=0.
    \end{cases}
\]
If $k(n)>0$, then $n\in E_{k(n)}$, so the midpoint defining $z_n$
lies strictly between $a_{k(n)}(n)$ and $b_{k(n)}(n)$.
Fix $j\ge1$. For every $n\in E_j$, we have $k(n)\ge j$, and hence
\[
    \ell_{j,n}\le a_{k(n)}(n)<z_n<b_{k(n)}(n)\le r_{j,n}.
\]
Thus $E_j\subseteq\{n:\ell_{j,n}<z_n<r_{j,n}\}$, and the latter set
belongs to $\mathcal U$. Consequently $z=[(z_n)]$ satisfies
$\ell_j<z<r_j$ for every $j$, proving $L<z<R$.

If $L\ne\varnothing$ and $R=\varnothing$, enumerate $L$ and choose
representatives as above, and set
\[
    z_n=1+\max_{1\le j\le n}\ell_{j,n}.
\]
For every fixed $j$, the inequality $\ell_{j,n}<z_n$ holds for all
$n\ge j$, so $\ell_j<[(z_n)]$. Similarly, if $L=\varnothing$ and
$R\ne\varnothing$, use
\[
    z_n=\min_{1\le j\le n}r_{j,n}-1
\]
to obtain $[(z_n)]<r_j$ for every $j$.
\end{proof}

We need the following property of cofinality of cuts in a totally ordered vector space with countable set of Archimedean classes. This follows from the cofinality formulas for cuts in divisible ordered abelian groups given in \cite[Theorems 4.5 and 4.8]{KuhlmannNart2022}, but we include a direct proof for completeness.

\begin{lema}\label{cpi:countable-cuts}
\leanref{FreeProducts/CountableSpectra/OrderedExtension.lean}{280}{countable_cut_points}
Let $V$ be a totally ordered real vector space such that $\Gamma_V$ is
at most countable. Let $D$ be a vector subspace of $V$ and
$b\in V\setminus D$. Then
\[
    L_b=\{d\in D:d<b\},
    \qquad
    R_b=\{d\in D:d>b\}
\]
have, respectively, a countable cofinal subset and a countable coinitial subset.
\end{lema}

\begin{proof}
We use the terminology for Archimedean classes introduced in Section \ref{subsec:countable spectrum}. The nonempty set
\[
    \Delta=\{v(b-d):d\in D\}\subseteq\Gamma_V
\]
is at most countable.

Suppose first that $\Delta$ has no largest element. Choose $d_n\in D$
such that
\[
    \delta_n=v(b-d_n)
\]
is strictly increasing and cofinal in $\Delta$. Set
\[
    h_n=|d_{n+1}-d_n|,
    \qquad
    \ell_n=d_n-2h_n,
    \qquad
    r_n=d_n+2h_n.
\]
Since $\delta_{n+1}>\delta_n$, the valuation identities give
\[
    v(d_{n+1}-d_n)=\delta_n,
    \qquad
    v(b-d_{n+1})>\delta_n.
\]
Consequently,
\[
    \ell_n<b<r_n,
    \qquad
    v(b-\ell_n)=v(r_n-b)=\delta_n.
\]
Indeed, write $b-d_n=(d_{n+1}-d_n)+(b-d_{n+1})$ and observe that
the second summand is infinitesimal relative to $h_n$. Both
$2h_n+(b-d_n)$ and $2h_n-(b-d_n)$ are therefore positive and have the
same valuation as $h_n$.

Given $d\in L_b$, choose $n$ with $\delta_n>v(b-d)$. Then
$0<b-\ell_n<b-d$, so $d<\ell_n$. Similarly, given $d\in R_b$, some $n$
satisfies $r_n<d$. Hence $\{\ell_n:n\ge1\}$ and $\{r_n:n\ge1\}$ have the
required properties.

Suppose now that $\Delta$ has largest element
$\delta=v(b-d_0)$, where $d_0\in D$. We claim that
\begin{equation}\label{cpi:new-class}
    \delta\notin v(D\setminus\{0\}).
\end{equation}
If this failed, choose $u\in D$, $u>0$, with $v(u)=\delta$, and put
$w=b-d_0$. Since $|w|\asymp u$, the number
\[
    a=\sup\{s\in\mathbb R:su\le w\}
\]
is a well-defined real number. For every $\eta>0$,
\[
    (a-\eta)u<w<(a+\eta)u,
\]
and thus $|w-au|<\eta u$. The vector $w-au$ is nonzero, since
$b\notin D$, and therefore
\[
    v\bigl(b-(d_0+au)\bigr)>\delta,
\]
contrary to the choice of $\delta$. This proves
\eqref{cpi:new-class}.

Replacing $b$ by $-b$ if necessary, it suffices to consider $w>0$.
For each $\gamma\in v(D\setminus\{0\})$, choose $u_\gamma\in D$ with
$u_\gamma>0$ and $v(u_\gamma)=\gamma$. Define
\[
\begin{split}
    L_0&=\{d_0\}\cup
       \{d_0+n u_\gamma:\gamma>\delta,\ n\ge1\},\\
    R_0&=\{d_0+n^{-1}u_\gamma:\gamma<\delta,\ n\ge1\}.
\end{split}
\]
These sets are at most countable and satisfy
$L_0\subseteq L_b$ and $R_0\subseteq R_b$.

To verify cofinality, let $d\in L_b$. If $d\le d_0$, use $d_0\in L_0$.
Otherwise $0<d-d_0<w$, so \eqref{cpi:new-class} implies
$v(d-d_0)>\delta$. A sufficiently large integer multiple of the chosen
representative of this class dominates $d-d_0$. Thus $d\le\ell$ for
some $\ell\in L_0$. Similarly, if $d\in R_b$, then
$d-d_0>w>0$ and $v(d-d_0)<\delta$. A sufficiently small positive
multiple $n^{-1}u_{v(d-d_0)}$ is less than $d-d_0$, proving that $R_0$
is coinitial in $R_b$.
\end{proof}

\begin{lema}\label{cpi:ordered-extension}
\leanref{FreeProducts/CountableSpectra/OrderedExtension.lean}{412}{ordered_lattice_embedding_extension}
Let $Y$ be a totally ordered real vector space with the countable
interpolation property. If $V$ is a
totally ordered real vector space with at most countably many
Archimedean classes, then every ordered linear embedding $j:D\to Y$,
defined on a vector subspace $D\subseteq V$, extends to an ordered
linear embedding $V\to Y$.
\end{lema}

\begin{proof}
By Zorn's lemma, choose a maximal extension $S:X\to Y$ of $j$ among
ordered linear embeddings defined on vector subspaces of $V$.
Suppose that $b\in V\setminus X$. By
Lemma~\ref{cpi:countable-cuts}, there are at most countable sets
$L_0\subseteq\{e\in X:e<b\}$ and
$R_0\subseteq\{e\in X:e>b\}$ that are cofinal and
coinitial in these sets, respectively. Countable interpolation yields $z\in Y$ with
\[
    S(L_0)<z<S(R_0).
\]
By cofinality and coinitiality,
\[
    S(e)<z\quad(e\in X,\ e<b),
    \qquad
    z<S(e)\quad(e\in X,\ b<e).
\]
Define
\[
    S'(e+tb)=S(e)+tz
    \qquad(e\in X,\ t\in\mathbb R).
\]
Since $b\notin X$, this defines a linear extension of $S$. If $t\ne0$,
the sign of $e+tb$ is determined by the comparison of $b$ with
$-e/t$. The displayed inequalities show that $S'(e+tb)$ has exactly
the same sign. For $t=0$, this follows from the corresponding property
of $S$. Hence $S'$ is an ordered linear embedding, contradicting
maximality. Therefore $X=V$.
\end{proof}

\begin{theo}\label{cpi:ultrapower}
\leanref{FreeProducts/CountableSpectra/Ultrapowers.lean}{177}{orderedUltrapower_isCountablePrimeInjective}
For every nonprincipal ultrafilter $\mathcal U$ on $\mathbb N$, the
totally ordered vector lattice $X_{\mathcal U}$ is
$\mathcal C$-injective. \leanref{FreeProducts/CountableSpectra/OrderedExtension.lean}{500}{isCountablePrimeInjective_of_countableInterpolation}
More generally, every totally ordered real
vector space with the countable interpolation property has this property.
\end{theo}
\begin{proof}
Let $X$ have the countable interpolation property, let $A\subset B\in\mathcal C$,
and let $T:A\to X$ be a nonzero lattice homomorphism.
Then $P=\ker T$ is a proper prime of $A$, because
$A/P$ is lattice isomorphic to the totally ordered space $T(A)$.
By Proposition~\ref{prop:prime-ideals-extcon}, choose a proper prime
$Q$ of $B$ with $Q\cap A=P$.

Let $q:B\to V=B/Q$ be the quotient map. By
Lemmas~\ref{cpi:hereditary} and~\ref{cpi:chain-cardinality}, the
totally ordered space $V$ has at most countably many Archimedean
classes. The map
\[
 j_0:q(A)\longrightarrow X,\qquad j_0(q(a))=T(a),
\]
is a well-defined ordered linear embedding. By
Lemma~\ref{cpi:ordered-extension}, it extends to an ordered linear
embedding $j:V\to X$. Such a map between totally ordered spaces is a
lattice homomorphism, so $\widetilde T=jq$ extends $T$ as required.
The case of $X_{\mathcal U}$ follows from Lemma~\ref{cpi:interpolation}.
\end{proof}

\begin{cor}\label{cpi:strictness}
\leanref{FreeProducts/Counterexamples/HyperrealDisposition.lean}{35}{hyperreal_injective_not_universalDisposition}
$\mathcal C$-injectivity does not imply universal disposition for
$\mathcal C$.
\end{cor}
\begin{proof}
The totally ordered lattice $X_{\mathcal U}$ cannot contain a lattice
copy of coordinatewise ordered $\mathbb R^2$, since it has no pair of
disjoint nonzero positive elements. But $\mathbb R^2$ has exactly two
proper primes, and universal disposition applied to
$0\subset\mathbb R^2$ requires such a copy.
\end{proof}

\begin{obs}
\leanref{FreeProducts/CountableSpectra/OrderedExtension.lean}{517}{exists_extension_of_countable_prime_quotient_classes}
The proof of Theorem~\ref{cpi:ultrapower} only requires each prime
quotient $B/Q$ to have countably many Archimedean classes. It imposes
no countability assumption on the Hamel dimension of $A$ or $B$.
The extension property is preserved by products and by retractions
that are lattice homomorphisms. \leanref{FreeProducts/CountableSpectra/InjectiveProperties.lean}{249}{isCountablePrimeInjective_pi}
For products, extend each coordinate
separately; \leanref{FreeProducts/CountableSpectra/InjectiveProperties.lean}{266}{isCountablePrimeInjective_of_retract}
for a retract, compose an extension with the retraction.
\end{obs}

\subsection{Properties of the $\mathcal C$-injective vector lattices}

\begin{prop}\label{cpi:local-properties}
Let $Z\neq0$ be $\mathcal C$-injective. Then:
\begin{enumerate}
\item \leanref{FreeProducts/CountableSpectra/InjectiveProperties.lean}{110}{exists_forall_smul_le_of_isCountablePrimeInjective}
For every $z>0$ there is
$y\in Z$ such that $nz\leq y$ for all $n\geq1$.
\item \leanref{FreeProducts/CountableSpectra/InjectiveProperties.lean}{123}{not_isVLArchimedean_of_isCountablePrimeInjective}
$Z$ is non-Archimedean.
\item \leanref{FreeProducts/CountableSpectra/InjectiveProperties.lean}{134}{not_strongOrderUnit_of_isCountablePrimeInjective}
$Z$ has no strong order unit.
\item \leanref{FreeProducts/CountableSpectra/InjectiveProperties.lean}{213}{exists_disjoint_positive_sum_of_universalDisposition}
If $U$ is of universal disposition for $\mathcal C$, then every
$u>0$ is the sum of two nonzero disjoint positive elements.
\item \leanref{FreeProducts/CountableSpectra/InjectiveProperties.lean}{227}{exists_positive_forall_nsmul_le_of_universalDisposition}
If $U$ is of universal disposition for $\mathcal C$, then for every
$u>0$ there is $x>0$ with $nx\leq u$ for all $n\geq1$.
\end{enumerate}
\end{prop}

\begin{proof}
The lexicographic plane $\mathbb R^2_{\mathrm{lex}}$, with its first
coordinate dominant, has two proper primes. Extend the homomorphism
$(0,t)\mapsto tz$ to this plane and let $y$ be the image of $(1,0)$.
Then $y\geq nz$ for all $n$. This forbids the Archimedean property. If $z$
were a strong order unit, $y\leq rz$ for some real $r$, also a
contradiction.

For the first assertion about $U$, extend $t\mapsto tu$ along the
diagonal embedding $t\mapsto(t,t)$ into coordinatewise
$\mathbb R^2$. The images of its coordinate vectors are the required
disjoint summands. For the second, extend $t\mapsto tu$ along
$t\mapsto(t,0)$ into $\mathbb R^2_{\mathrm{lex}}$. The image of
$(0,1)$ is the required nonzero element.
\end{proof}

\begin{prop}\label{cpi:chain-universality}
\leanref{FreeProducts/CountableSpectra/InjectiveProperties.lean}{385}{exists_injective_hom_of_isCountablePrimeInjective}
Every nonzero $\mathcal C$-injective vector lattice $Z$ contains a
lattice copy of every totally ordered member of $\mathcal C$.
\leanref{FreeProducts/CountableSpectra/Basic.lean}{529}{continuum_le_cardinal_of_nontrivial}
In particular, $|Z|\geq\mathfrak c$. \leanref{FreeProducts/CountableSpectra/Ultrapowers.lean}{167}{orderedUltrapower_cardinal}
This lower bound is attained by
$X_{\mathcal U}$.
\end{prop}

\begin{proof}
Let $C\in\mathcal C$ be totally ordered and form
$W=C\oplus_{\mathrm{lex}}\mathbb R$, with the $C$ coordinate dominant.
Its proper ideals are $0$ and $J\oplus\mathbb R$, where $J$ ranges
over the proper ideals of $C$. Thus $W\in\mathcal C$, and every
nonzero ideal of $W$ contains $e=(0,1)$. For $z>0$ in $Z$, extend
$te\mapsto tz$ to $W$. The extension has zero kernel, because a
nonzero kernel would contain $e$. Restriction to $C\times\{0\}$ gives
the desired embedding. Taking the lexicographic product from
Example~\ref{cpi:lex-product} yields $|Z|\geq\mathfrak c$; the
ultrapowers $X_{\mathcal U}$ have cardinality exactly $\mathfrak c$.
\end{proof}

For an ordinal $\alpha>0$, let $H_\alpha=\mathbb R^{(\alpha)}$,
ordered by the sign of the coefficient at the largest support index.
Its proper ideals are precisely the initial-coordinate spaces
$H_\beta$, $\beta<\alpha$: a vector with largest support index
$\gamma$ generates the ideal $H_{\gamma+1}$, and all ideals are
unions of such ideals. In particular,
$|\operatorname{Spec}_p(H_\alpha)|=|\alpha|$, and every nonzero ideal
of $H_\alpha$ contains $e_0$.

\begin{theo}\label{cpi:large-spectrum}
\leanref{FreeProducts/CountableSpectra/InjectiveProperties.lean}{588}{exists_omega_one_colex_embedding}
Every nonzero $\mathcal C$-injective vector lattice $Z$ contains
a lattice copy of $H_{\omega_1}$.
\leanref{FreeProducts/CountableSpectra/InjectiveProperties.lean}{24}{continuum_le_spectrum_of_isCountablePrimeInjective}
Its prime spectrum satisfies
\[
 |\operatorname{Spec}_p(Z)|\geq\mathfrak c.
\]
\leanref{FreeProducts/CountableSpectra/InjectiveProperties.lean}{36}{countable_spectrum_not_injective}
In particular, no nonzero member of $\mathcal C$ is
$\mathcal C$-injective.
\end{theo}
\begin{proof}
Fix $z>0$ in $Z$ and embed $H_1$ by $e_0\mapsto z$. Construct
compatible embeddings $T_\alpha:H_\alpha\to Z$ for
$1\leq\alpha<\omega_1$. At a successor ordinal $H_{\alpha+1}$ has
countable spectrum, so relative injectivity extends $T_\alpha$ to a
homomorphism on $H_{\alpha+1}$. This extension is injective: any
nonzero kernel would contain $e_0$. At limit ordinals take unions.
The final union embeds $H_{\omega_1}$ into $Z$.

For the stronger spectrum bound, let $F_n$ consist of $0,1$ and the
first $n$ terms of an enumeration of $\mathbb Q\cap(0,1)$. Give
$H(F_n)=\mathbb R^{F_n}$ the order determined by its largest nonzero
coordinate. Successively extending $e_0\mapsto z$ gives compatible
embeddings $H(F_n)\to Z$: these test lattices have finite spectra,
and a nonzero kernel would again contain $e_0$. Their union embeds
$H(\mathbb Q\cap[0,1])=\mathbb R^{(\mathbb Q\cap[0,1])}$ into $Z$.
For each $t\in(0,1)$, the span of the $e_q$ with $q<t$ is a proper
prime ideal of this chain. These ideals are distinct, so
Lemma~\ref{cpi:hereditary} gives
$|\operatorname{Spec}_p(Z)|\geq\mathfrak c$.
\end{proof}

\subsection{An obstruction using Archimedean vector lattices}

Proposition~\ref{cpi:local-properties} uses a non-Archimedean vector
lattice. The obstruction persists when all vector lattices involved are required
to be Archimedean and to have countable spectrum.

Let $\mathbf1=(1,1,\ldots)$, $b=(1,2,3,\ldots)$, and let $c_{00}$
denote the finitely supported real sequences. Give
\[
 A=c_{00}+\mathbb R\mathbf1,\qquad
 B=c_{00}+\mathbb R\mathbf1+\mathbb R b
\]
their coordinatewise orders. These are the eventually constant and
eventually affine sequences, respectively. Since an affine sequence
has constant sign eventually, both spaces are vector sublattices of
$\mathbb R^{\mathbb N}$ and hence are Archimedean.

\begin{lema}\label{cpi:sequence-primes}
Writing $P_n^A=\{a\in A:a_n=0\}$ and
$P_n^B=\{x\in B:x_n=0\}$, we have:
\begin{enumerate}
\item \leanref{FreeProducts/Counterexamples/SequenceLattices.lean}{521}{eventuallyConstant_spectrum}
\[
 \operatorname{Spec}_p(A)=\{P_n^A:n\geq1\}\cup\{c_{00}\}.
\]
\item \leanref{FreeProducts/Counterexamples/SequenceLattices.lean}{503}{eventuallyAffine_spectrum}
\[
 \operatorname{Spec}_p(B)=\{P_n^B:n\geq1\}\cup\{c_{00},A\}.
\]
\item \leanref{FreeProducts/Counterexamples/SequenceLattices.lean}{375}{sequence_lattices_countable_spectrum}
$A,B\in\mathcal C$.
\item \leanref{FreeProducts/Counterexamples/SequenceLattices.lean}{348}{affine_sublattice_eq_top}
$B$ is generated as a
vector lattice by $\mathbf1$ and $b$.
\end{enumerate}
\end{lema}

\begin{proof}
In either space, if a prime ideal $P$ does not contain $e_n$, disjointness with $e_n$
forces every vector whose $n$th coordinate is zero to belong to $P$. Conversely,
if $x\in P$, the inequality $|x_n|e_n\leq|x|$ forces $x_n=0$.
Thus $P$ is the corresponding coordinate kernel. All other primes
contain $c_{00}$.

Now $A/c_{00}\cong\mathbb R$, while
$B/c_{00}\cong\mathbb R^2_{\mathrm{lex}}$ under
$[a\mathbf1+s b]\mapsto(s,a)$. The latter has exactly the two proper
ideals $0$ and $\{0\}\times\mathbb R$, whose inverse images are
$c_{00}$ and $A$. This gives the stated descriptions of the spectra. Finally,
\[
 e_n=(\mathbf1-|b-n\mathbf1|)^+\qquad(n\geq1),
\]
so $\mathbf1,b$ generate $B$.
\end{proof}

\begin{theo}\label{cpi:archimedean-obstruction}
\leanref{FreeProducts/Counterexamples/NonextendibleHom.lean}{45}{archimedean_nonextendible_hom}
No nonzero Archimedean vector lattice is injective relative to
$\mathcal C\cap\mathcal{AVL}$. \leanref{FreeProducts/Counterexamples/ArchimedeanDisposition.lean}{84}{archimedean_not_universalDisposition}
There is no Archimedean vector lattice
of universal disposition for this class either.
\end{theo}
\begin{proof}
If $Z\neq0$ is Archimedean, choose $u>0$ in $Z$ and define the lattice
homomorphism
\[
 T:A\to Z,\qquad T(a)=(\lim_{n\to\infty}a_n)u.
\]
For every integer $m\geq1$, the vector
$a^{(m)}=b\wedge m\mathbf1$ lies in $A$ and has eventual value $m$.
Any extension $\widetilde T:B\to Z$ would satisfy
\[
 0\leq mu=T(a^{(m)})\leq\widetilde T(b)\qquad(m\geq1),
\]
contrary to the Archimedean property.

For universal disposition, put $D=A\oplus\mathbb R$ with
coordinatewise order, and consider the embeddings
\[
 i:A\to B,\quad i(a)=a,
 \qquad j:A\to D,\quad j(a)=(a,\lim a).
\]
The lattice $D$ is Archimedean and has countable spectrum: the proper
primes of a finite direct product come from a prime in one factor.
If homomorphisms $f:B\to Z$ and $g:D\to Z$ into an Archimedean
lattice satisfy $fi=gj$, then
\[
 f(b)\geq f(a^{(m)})=g(a^{(m)},m)\geq m g(0,1)\geq0
 \qquad(m\geq1).
\]
Hence $g(0,1)=0$, and $g$ cannot be injective. If an Archimedean $U$
had the stated universal disposition, first embed $D$ into $U$ using
$0\subset D$ and then extend the resulting embedding of $A$ along
$i:A\to B$. This would give such an amalgam with $g$ injective, a
contradiction.
\end{proof}

\leanref{FreeProducts/Counterexamples/CountableSpectrumPushout.lean}{68}{exists_countable_spectrum_nonarchimedean_pushout}
In the ordinary $\mathcal{VL}$ push-out of $i$ and $j$,
Proposition~\ref{prop:pushout-vl} makes the image $e$ of $(0,1)\in D$
nonzero. If $y$ is the image of $b\in B$, the same calculation gives
$0<me\leq y$ for every $m\geq1$. Thus even a diagram of Archimedean
lattices with countable spectra can force a non-Archimedean push-out;
its Archimedeanization sends $e$ to zero.

The existence results above therefore belong to $\mathcal{VL}$.
They provide no uniqueness theorem for universal disposition relative
to $\mathcal C$: a back-and-forth construction would require control
of the intermediate sublattices, and
Example~\ref{cpi:closure-failures} shows why countable spectrum alone
does not supply the usual closure assumptions.

\section*{Acknowledgements}
Research partially supported by grants PID2024-162214NB-I00 and CEX2023-001347-S funded by MCIN/AEI/10.13039/501100011033. Martínez-Fernández was also supported by grant PREX2023-000080 funded by MCIN/AEI/10.13039/501100011033.



\bibliographystyle{amsplain-initials}
\bibliography{mainbib}

\end{document}